\documentclass[reqno,11pt,letterpaper]{amsart}

\usepackage[dvipsnames]{xcolor}
\usepackage{xcolor}
\usepackage{amsthm}
\usepackage{amsmath}
\usepackage{amssymb}
\usepackage{lmodern}
\usepackage{mathrsfs}
\usepackage{amsfonts}
\usepackage{enumerate}
\usepackage{txfonts}
\usepackage{bm}

\usepackage[numbers]{natbib}
\usepackage{ulem}

\usepackage{graphicx}

\numberwithin{equation}{section}

\usepackage[colorlinks=true]{hyperref}

\def\ba{\bar{a}}

\def\bU{\bar{U}}

\def\bS{\bar{S}}
\def\bB{\bar{B}}
\def\bM{\bar{M}}

\def\brho{\bar{\rho}}

\def\bka{\bar{\kappa}}

\def\hV{\hat{V}}

\def\be{\begin{eqnarray}}
\def\ee{\end{eqnarray}}
\def\ba{\begin{aligned}}
\def\ea{\end{aligned}}
\def\bay{\begin{array}}
\def\eay{\end{array}}
\def\bca{\begin{cases}}
\def\eca{\end{cases}}
\def\p{\partial}
\def\no{\nonumber}

\def\f{\frac}
\def\th{\theta}

\def\lab{\label}
\def\b{\bigg}

\def\ka{\kappa}

\def\ga{\gamma}

\def\q{\quad}

\def\u{\mathbf{u}}
\def\h{\mathbf{h}}
\def\mlS{\mathcal{S}}

\newtheorem{theorem}{Theorem}[section]
\newtheorem{lemma}[theorem]{Lemma}

\newtheorem{proposition}[theorem]{Proposition}

\theoremstyle{definition}

\theoremstyle{remark}

\allowdisplaybreaks[1]

\begin{document}

\title[super-Alfv\'enic MHD shocks]{On three dimensional steady super-Alfv\'{e}nic magnetohydrodynamics shocks with aligned fields}

\author[S.K. Weng]{Shangkun Weng}
\address[Shangkun Weng]{\newline School of Mathematics and Statistics, Wuhan University, Wuhan, Hubei Province, 430072, People's Republic of China.}
\email{skweng@whu.edu.cn}

\author[W.G. Yang]{Wengang Yang}
\address[Wengang Yang]{\newline School of Mathematics and Statistics, Wuhan University, Wuhan, Hubei Province, 430072, People's Republic of China.}
\email{yangwg@whu.edu.cn}

\begin{abstract}
The coupled motion between the hydrodynamic flow and magnetic field introduces significant complexity into the structure of the magnetohydrodynamic (MHD) equations. A key factor contributing to this complexity is the presence of Alfv\'en waves, which critically influences the character of the flow and makes the problem considerably more challenging. Within the framework where the magnetic field is everywhere parallel to the flow velocity, we give an effective decomposition of the steady MHD equations in terms of the deformation tensor and the modified vorticity, where the modification in the vorticity is to record the effect of the Lorentz force on the velocity field. The existence and structural stability of the super-Alfv\'{e}nic cylindrical transonic shock solutions for the steady MHD equations are established under three-dimensional perturbations of the incoming flow and the exit total pressure (kinetic plus magnetic).

\end{abstract}

\keywords{MHD shock, elliptic-hyperbolic mixed, sup-Alfv\'{e}nic, Rankine-Hugoniot conditions.}
\subjclass[2020]{35L67, 35M12, 76N10, 76N15}
\date{}
\maketitle

\section{Introduction and main results}\label{sec1}
The dynamics of a compressible and inviscid magnetohydrodynamics (MHD) fluid are described by the following equations:
\begin{equation}\label{mhd_t0}
\begin{cases}
\p_t\rho+\nabla\cdot(\rho \u)=0,\\
\p_t(\rho \u)+\nabla\cdot(\rho \u\otimes \u-\h\otimes \h)+\nabla(P+\frac{1}{2}|\h|^2)=0,\\
\p_t(\rho (\frac{1}{2}|\u|^2+e)+\frac{1}{2}|\h|^2)+\nabla\cdot(\rho \u(\frac{1}{2}|\u|^2+e+\frac{P}{\rho})+\h\times(\u\times \h))=0,\\
\p_t\h-\nabla\times(\u\times\h)=0,
\end{cases}
\end{equation}
and 
\begin{equation}\label{mhd_t1}
\nabla\cdot\h=0,
\end{equation}
where ${\bf u}$, $\rho$, $P$, $e$, and ${\bf h}$ represent the velocity, density, pressure, internal energy, and the magnetic field, respectively. The system of MHD combines principles from fluid mechanics and electromagnetism to form a unified theory for studying electrically conducting fluids. The MHD equations \eqref{mhd_t0}-\eqref{mhd_t1} are applicable to a broad range of physical regimes, including plasmas, astrophysics, and controlled nuclear fusion \cite{Davidson2017}.
Mathematically, the analysis of fundamental nonlinear waves, such as shocks, rarefaction waves, and contact discontinuities, constitutes a central theme in the study of multidimensional hyperbolic conservation laws. In the MHD framework \cite{BT2002}, significant research has examined the existence, uniqueness, and stability of fundamental waves, focusing particularly on shocks \cite{Trak2003CMP}, vortex sheets \cite{ChenW2008,Tark2009ARMA}, and contact discontinuities \cite{MTT2018,TrakWang2022,WangXin2024}. 

In this paper, we investigate stationary shocks in a de Laval nozzle composed of a concentric cylindrical segment. The flow is governed by the three-dimensional steady compressible MHD equations:
\be\label{mhd}\begin{cases}
\text{div }(\rho {\bf u})=0,\\
\text{div }(\rho {\bf u}\otimes {\bf u}+ P I_3) =\text{curl }{\bf h}\times {\bf h},\\
\text{div }(\rho\u B+\h\times({\bf u}\times {\bf h})) =0 ,\\
\text{curl }({\bf u}\times {\bf h})=0,\\
\text{div }{\bf h}=0,
\end{cases}\ee
where  $B=\frac{|{\bf u}|^2}{2}+e+\frac{P}{\rho}$ denotes the Bernoulli quantity. For simplicity, we consider the polytropic gas, whose equation of state and the internal energy are
\begin{equation}\no
P=S\rho^{\gamma}\quad \text{and}\quad
e=\frac{P}{(\gamma-1)\rho},\ \ \ \gamma>1,
\end{equation}
respectively, where $\gamma> 1$ and $S$ is the entropy.

A discontinuity front is a surface at which some or all of the above quantities have jump discontinuities. Let the discontinuity front be given by $\mlS(x_1,x_2,x_3)=0$. The two states on the two sides of the front are connected by the Rankine-Hugoniot jump conditions, which are obtained by integrating the equations \eqref{mhd} across the front and are as follows:
\be\no\begin{cases}
\nabla \mlS\cdot [\rho {\bf u}]=0,\\
\nabla \mlS\cdot [\rho {\bf u}\times {\bf u}-{\bf h}\times {\bf h}+ (P+\frac12 |{\bf h}|^2)I_{3\times 3}]=0,\\
\nabla \mlS\cdot [\rho(e+\frac12 |{\bf u}|^2){\bf u}+P {\bf u}-({\bf u}\times {\bf h})\times {\bf h}]=0,\\
\nabla \mlS\times [{\bf u}\times {\bf h}]=0,\\
\nabla \mlS\cdot [{\bf h}]=0.\end{cases}
\ee

The study of transonic shocks for the Euler equations can be traced back to the seminal work \cite{CF1948} of Courant and Friedrichs in 1948. Compared to the steady Euler equations, the Lorentz force induced by the magnetic field in MHD flow introduces a fundamental difference from pure gas dynamics by facilitating the anisotropic propagation of small disturbances. A detailed description of the MHD analogues of shocks and sound waves was exhibited by De Hoffmann and Teller, Friedrichs \cite{DT1950,Fried1954}, including the Alfv\'{e}n wave and simple waves. Bazer and Ericson \cite{BE1959} investigated the one-dimensional steady motion of \eqref{mhd} and gave a complete classification of all physically admissible solutions of the MHD discontinuity relations. Define $V=\frac{1}{\rho},i_n=\rho {\bf u}\cdot {\bf n}$, $h_n={\bf h}\cdot {\bf n}$ and ${\bf h}_{\tau}={\bf h}-h_n {\bf n}$, where ${\bf n}$ is the outer normal to the discontinuity front. The complete list of admissible solutions to the discontinuity relations is as follows (see \cite{FT2013}):
\begin{enumerate}[(1)]
\item $i_n=0, h_n\neq 0$ (contact discontinuity); here $[{\bf u}]=[{\bf h}]=[P]=0$. Other parameters, the density, the temperature can have arbitrary jumps.
\item $i_n=0, h_n=0$ (tangential discontinuity). The tangential component of the velocity and magnetic field can have arbitrary jumps. Thermodynamic parameters of the fluid can also display jumps with a constraint $[P+\frac{1}{2}|{\bf h}_{\tau}|^2]=0$.
\item $i_n\neq 0, [V]=0$ (Alfv\'{e}n or rotational discontinuity). All the internal energy, entropy, kinetic pressure and the magnitude of the magnetic field are continuous through the Alfv\'{e}n discontinuity. However, $[{\bf h}_{\tau}]\neq 0$, which means that the magnetic field ${\bf h}$ rotates by an arbitrary angle at the Alfv\'{e}n discontinuity.
\item $i_n\neq 0, [V]\neq 0$ (shock front). Then the Rankine-Hugoniot shock adiabat equation holds
\be\no
e_1-e_2+ \frac{1}{2}(p_1+p_2)(V_1-V_2)+\frac{1}{4}(V_1-V_2)(|{\bf h}_{\tau1}|-|{\bf h}_{\tau 2}|)^2=0,
\ee
where the index $2$ marks the downstream values, while $1$ marks the upstream values. As required by the second law of thermodynamics, the entropy, the pressure and the density increase:
\be\no 
s_2>s_1,\ \ P_2>P_1, \ \ \rho_2>\rho_1.
\ee
\end{enumerate}

For the stability analysis of transonic shocks in finitely long nozzles, two types of transonic shock solutions commonly serve as fundamental reference flows. The first type consists of two constant states with an arbitrarily located shock. The existence, uniqueness, and structural stability in nozzles under various boundary conditions were studied in \cite{ChenF2003,Chen2008,ChenYuan2008,XinYin2005,XinYin2008PJM} for multidimensional steady potential flow. Fang and Xin \cite{FangX2021} developed an elaborate approach to uniquely determine the position of the shock front for two-dimensional steady Euler equations in an almost flat nozzle.




The second type involves symmetric transonic shocks in divergent nozzles, such as radial shocks in angular sectors or spherical shocks in cones, where the shock position is uniquely determined by the exit pressure. The existence and stability of transonic shock solutions in divergent nozzles under general perturbations of the wall and exit pressure were established in \cite{LXY2009CMP,LXY2013}. Corresponding results for axisymmetric perturbations were subsequently examined in \cite{LXY2010JDE,WengXieXin2021}. In \cite{LiuXY2016}, the stability of spherically symmetric subsonic flows and transonic shocks in a spherical shell was established under certain “Structural Conditions” imposed on the background transonic shock solutions. Recently, the authors removed these structural conditions and proved the existence and stability of both cylindrical \cite{WengX24ArxicCyl} and spherical \cite{Weng25ArxivSph} transonic shocks under three-dimensional perturbations of the incoming flow and exit pressure. This was achieved by employing the deformation-curl decomposition \cite{WengX2019}, a refined reformulation of the Rankine-Hugoniot conditions, and the introduction of ``spherical projection coordinates". The authors in \cite{WengY2024,WengZZ2025} proved the structural stability of a transonic shock in a two-dimensional flat nozzle under an external force, revealing how a well-chosen force can exert a stabilizing effect on the shock. 

Clearly, the character of the steady Euler equations is fully determined by the Mach number, defined as $M^2 = |\mathbf{u}|^2 / c^2(\rho, S)$, where $c(\rho, S) = \sqrt{\partial_\rho P(\rho, S)}$ is the local sound speed. In contrast, for steady MHD flows, the type of the governing differential equations depends not only on the Mach number but also on another key dimensionless parameter: the Alfv\'{e}n number $A^2 = |\mathbf{u}|^2 / c_a^2$, where $c_a = \sqrt{|\mathbf{h}|^2 / \rho}$ denotes the Alfv\'{e}n wave speed. For an infinitely conducting accelerating transonic gas with the magnetic field parallel to the velocity everywhere, the gas must cross three transitions \cite{Chu1962,Grad1960} at $A^2+M^2=1, A^2=1$ and $M^2=1$, respectively. The steady MHD equations are elliptic-hyperbolic mixed (purely hyperbolic) if $(A^2-1)(M^2-1)(A^2+M^2-1)<0 (>0, \text{respectively})$. Consequently, the mathematical analysis of the steady compressible MHD equations is significantly more complex and challenging than that of the steady compressible Euler equations. 


As a preliminary investigation of the steady MHD shock, we focus on the case of aligned magnetic and velocity fields. Specifically, we assume that the magnetic field $\h$ and the velocity field $\u$ are everywhere parallel:
\be\label{ka}
{\bf h}=\kappa \rho {\bf u},
\ee
where $\kappa$ is a scalar function, then \eqref{mhd} simplifies to
\be\label{mhd1}\begin{cases}
\text{div }(\rho {\bf u})=0,\\
\text{div }(\rho {\bf u}\otimes {\bf u}+ P I_3) =\kappa \text{curl }(\kappa \rho{\bf u})\times (\rho{\bf u}),\\
\text{div }(\rho (\frac12|{\bf u}|^2 +e) {\bf u}+ P{\bf u}) =0,\\
\rho {\bf u}\cdot\nabla \kappa=0.
\end{cases}\ee

For our purpose, we introduce the cylindrical coordinates
\be\no
x_1=r\cos\theta,\ x_2=r\sin\theta,\ x_3=x_3,
\ee
and represent the velocity field as ${\bf u}(x)=U_1{\bf e}_r+U_2{\bf e}_\th+U_3{\bf e}_3$, where
\be\no
{\bf e}_r=(\cos\theta, \sin\theta, 0)^t, \ \ {\bf e}_\theta=(-\sin\theta,\cos\theta,0)^t,\ \ {\bf e}_3=(0,0,1)^t.
\ee
Then \eqref{mhd1} takes the following form in cylindrical coordinates
\begin{eqnarray}\label{mhd-cyl}
\begin{cases}
\p_r(\rho U_1)+\frac1r \rho U_1+\frac{1}{r}\p_{\theta} (\rho U_2) + \p_{x_3} (\rho U_3)=0,\\
(U_1\p_r +\frac{U_2}{r}\p_{\theta}+U_3\p_{x_3}) U_1+\frac1\rho \p_r P-\frac{U_2^2}r\\
\quad=-\kappa U_2\{(\p_r+\frac1{r})(\kappa \rho U_2)-\frac1{r}\p_{\theta}(\kappa \rho U_1)\}+\kappa U_3 \{\p_3(\kappa \rho U_1)-\p_r(\kappa \rho U_3)\},\\
(U_1\p_r +\frac{U_2}{r}\p_{\theta}+U_3\p_{x_3}) U_2+\frac1{r\rho} \p_{\theta} P+\frac{U_1U_2}r=\\
\quad=\kappa U_1\{(\p_r+\frac1{r})(\kappa \rho U_2)-\frac1{r}\p_{\theta}(\kappa \rho U_1)\}-\kappa U_3 \{\frac1{r}\p_{\theta}(\kappa \rho U_3)-\p_{x_3}(\kappa \rho U_2)\},\\
(U_1\p_r +\frac{U_2}{r}\p_{\theta}+U_3\p_{x_3}) U_3+\frac1\rho \p_{x_3}P=\\
\quad=\kappa U_2\{\frac{1}{r}\p_{\theta}(\kappa \rho U_3)-\p_{x_3}(\kappa \rho U_2)\}-\kappa U_1\{\p_3(\kappa \rho U_1)-\p_r(\kappa \rho U_3)\},\\
(U_1\p_r +\frac{U_2}{r}\p_{\theta}+U_3\p_{x_3}) B=0,\\
(U_1\p_r +\frac{U_2}{r}\p_{\theta}+U_3\p_{x_3}) \kappa=0.
\end{cases}
\end{eqnarray}
The flow region is assumed to be a part of a concentric cylinder described as
\be\no
\Omega=\{(r,\theta,x_3):r_1<r<r_2,(\theta,x_3)\in E\}, \ \ E:=(-\theta_0,\theta_0)\times (-1,1),
\ee
where $0<r_1<r_2<\infty,\theta_0\in (0,\frac{\pi}{2})$ are fixed positive constants.

Now we construct a class of cylindrically symmetric shock solutions with only nontrivial radial velocity to \eqref{mhd-cyl}, then \eqref{mhd-cyl} further reduces to
\begin{eqnarray}\label{mhd-cyl1}
\begin{cases}
(\bar{\rho} \bar{U})'(r)+\frac1r \bar{\rho} \bar{U}=0,\\
\bar{\rho} \bar{U}\bar{U}' +\bar{P}'(r)=0,\\
\bar{\rho} \bar{U} \bar{B}'(r)=0,\\
\bar{\rho} \bar{U}\bar{\kappa}'(r)=0.
\end{cases}
\end{eqnarray}
The corresponding Rankine-Hugoniot conditions and the physical entropy condition at the shock $r=r_s$ are
\be\label{cylsym-rh}\begin{cases}
[\bar{\rho} \bar{U}](r_s)=[\bar{\rho} \bar{U}^2+\bar{P}](r_s)=0,\\
[\bar{B}](r_s)=[\bar{\kappa}](r_s)=0,\quad \bar{S}_+>\bar{S}_-,
\end{cases}\ee
where $[f](r_s):=f(r_s+)-f(r_s-)$ denotes the jump of $f$ at $r=r_s$.

\begin{proposition}\label{background}
{\it Given the incoming supersonic flow $(\bar{U}_-(r_1){\bf e}_r, \bar{\rho}_-(r_1),\bar{S}_-,\bar\kappa)$ at $r=r_1$, where $\bar{U}_-(r_1)>0,\bar{\rho}_-(r_1)>0, \bar{S}_->0$ and $\bar{U}_-^2(r_1)>c^2(\bar{\rho}_-(r_1),\bar{S}_-)$. Then there exist two positive constants $P_1$ and $P_2$ depending only on the incoming supersonic flow and $r_1,r_2$, such that when the exit pressure $P_e\in (P_1, P_2)$, there exists a unique cylindrically symmetric shock solution
\be\no\begin{cases}
({\bar {\bf u}}_-,\bar{\rho}_-,\bar{S}_-,\bar{\kappa}_-):=(\bar{U}_-(r){\bf e}_r, \bar{\rho}_-(r),\bar{S}_-,\bka),\ \ &\text{in }(r_1,r_s),\\
({\bar {\bf u}}_+,\bar{\rho}_+,\bar{S}_+,\bar{\kappa}_+):=(\bar{U}_+(r){\bf e}_r, \bar{\rho}_+(r),\bar{S}_+,\bka),\ \ &\text{in }(r_s,r_2),
\end{cases}\ee
to \eqref{mhd-cyl1}, which satisfies the incoming supersonic flow and the exit pressure
\be\no
\bar{P}_+(r_2)=P_e
\ee
with a shock front at $r=r_s\in (r_1,r_2)$ satisfying \eqref{cylsym-rh}.  }
\end{proposition}

Later on, this special solution, ${\bf \overline{\Psi}}$, will be called the background solution. Clearly, one can extend the supersonic and subsonic parts of ${\bf \overline{\Psi}}$ in a natural way, respectively. With an abuse of notations, we still call the extended subsonic and supersonic solutions ${\bf\overline{\Psi}}_+$ and ${\bf\overline{\Psi}}_-$, respectively.  For detailed properties of this cylindrically symmetric transonic shock solution, we refer to \cite[Section 147]{CF1948} or \cite[Theorem 1.1]{XinYin2008JDE}. The main goal of this paper is to establish the structural stability of this cylindrically symmetric transonic shock solution under generic three-dimensional perturbations of suitable boundary conditions at the entrance and exit.

As we have discussed above, the presence of the magnetic field may greatly change the type of differential equations for the supersonic and subsonic flows. Here, we first consider the sup-Alfv\'{e}nic case, meaning that
\be\label{supA1}
&&\bar{A}_-^2(r)=\frac{1}{\bar{\kappa}^2 \bar{\rho}_-(r)}>1, \forall r\in [r_1,r_s],\\\label{supA2}
&&\bar{A}_+^2(r)=\frac{1}{\bar{\kappa}^2 \bar{\rho}_+(r)}>1, \forall r\in [r_s,r_2].
\ee
This condition is equivalent to
\be\label{supA0}
\bar{\kappa}^2<\min\big\{\min\limits_{r\in [r_1,r_s]}(\bar{\rho}_-(r))^{-1},\min\limits_{r\in [r_s,r_2]}(\bar{\rho}_+(r))^{-1}\big\}.
\ee

We will verify later that if $\bar{\kappa}^2$ satisfies \eqref{supA1}, then \eqref{mhd-cyl} for the upstream supersonic flows is purely hyperbolic. If $\bar{\kappa}^2$ satisfies \eqref{supA2}, then \eqref{mhd-cyl} for the downstream subsonic flows is elliptic-hyperbolic mixed. Under these assumptions for $\bar{\kappa}$, we now formulate suitable boundary conditions at the entrance and exit of the nozzle to find shock solutions to \eqref{mhd-cyl} which are close to the background cylindrical symmetric shock solutions.

Let the incoming supersonic flow at the inlet $r=r_1$ be prescribed as
\begin{equation}\lab{super1}
{\bf \Psi}_-(r_1, \theta,x_3)={\bf \overline{\Psi}}_-(r_1)+ \epsilon (U_{10},U_{20},U_{30}, P_{0},S_{0},\kappa_{0})(\theta,x_3),
\end{equation}
where
$(U_{10},U_{20},U_{30}, P_{0},S_{0},\kappa_{0})\in (C^{2,\alpha}(\overline{E}))^6$.
The flow satisfies the slip condition ${\bf u}\cdot {\bf n}$=0 on the nozzle wall, where ${\bf n}$ is the outer normal of the nozzle wall, which in the cylindrical coordinates, can be written as
\be\lab{slip1}\begin{cases}
U_2(r,\pm \theta_0,x_3) =0\quad \ &\forall (r, x_3)\in [r_1,r_2]\times [-1, 1],\\
U_3(r,\theta,\pm 1) =0\quad \ &\forall (r,\theta,)\in [r_1,r_2]\times[-\theta_0, \theta_0].
\end{cases}\ee

At the exit of the nozzle $\Gamma_o:=\{(r_2,\theta,x_3): (\theta,x_3)\in E\}$, different from the pure gas dynamics case, we should prescribe the total pressure
\be\lab{pressure}
(P+\frac{1}{2}|{\bf h}|^2)(r_2,\theta,x_3)= (\bar{P}+\frac12\bar{\kappa}^2\bar{\rho}^2\bar{U}^2)(r_2) + \epsilon T_{e}(\theta,x_3),
\ee
here $T_{e}\in C^{2,\alpha}(\overline{E})$ satisfies the compatibility conditions
\be\label{pressure-cp}\begin{cases}
\p_{\theta} T_{e}(\pm \theta_0,x_3)=0,\ \  \ &\forall x_3\in [-1,1],\\
\p_{x_3} T_{e}(\theta,\pm 1)=0,\ \  \ &\forall \theta\in [-\theta_0,\theta_0].
\end{cases}\ee

The problem is to find a piecewise smooth solution $\bm \Psi$ to \eqref{mhd-cyl} supplemented with the boundary conditions \eqref{super1}, \eqref{slip1}, and \eqref{pressure}, which jumps only at a shock front $\mlS: r=\xi(\theta,x_3), (\theta,x_3) \in E$. More precisely, we would construct functions
\begin{equation}\nonumber
\begin{cases}
(U_{1-}, U_{2-}, U_{3-}, P_-, S_-, \kappa_-),\quad\text{in}\,\, \Omega_-=\{r_1< r<\xi(\theta,x_3),(\theta,x_3)\in E\},\\
(U_{1+}, U_{2+}, U_{3+}, P_+, S_+, \kappa_+),\quad\text{in}\,\, \Omega_+=\{\xi(\theta,x_3)<r<r_2,(\theta,x_3)\in E\},
\end{cases}
\end{equation}
solve the equations \eqref{mhd-cyl} in $\Omega_{\pm}$ and satisfy the Rankine-Hugoniot conditions on the shock front $r=\xi(\theta, x_3)$:
\begin{eqnarray}\label{rh}
\begin{cases}
[\rho U_1]-\frac{1}{\xi}\p_{\theta}\xi [\rho U_2]-\p_{x_3}\xi[ \rho U_3]=0,\\
[\rho U_1^2+P+\frac{1}{2}\kappa^2 \rho^2(U_2^2+U_3^2-U_1^2)]-\frac{1}{\xi}\p_{\theta}\xi[(1-\ka^2\rho)\rho U_1U_2]\\
\q\q\q-\p_{x_3}\xi[(1-\ka^2\rho)\rho U_1U_3]=0,\\
[(1-\ka^2\rho)\rho U_1U_2]-\frac{1}{\xi}\p_{\theta}\xi[\rho U_2^2+P+\frac{1}{2}\kappa^2 \rho^2(U_1^2+U_3^2-U_2^2)]\\
\q\q\q-\p_{x_3}\xi[(1-\ka^2\rho)\rho U_2U_3]=0,\\
[(1-\ka^2\rho)\rho U_1U_3]-\frac{1}{\xi}\p_{\theta}\xi[(1-\ka^2\rho)\rho U_2U_3]\\
\q\q\q-\p_{x_3}\xi[ \rho U_3^2+P+\frac{1}{2}\kappa^2 \rho^2(U_1^2+U_2^2-U_3^2)]=0,\\
[B]=[\kappa]=0,
\end{cases}
\end{eqnarray}
and the physical entropy condition
\be\label{entropy}
S_+(\xi(\theta, x_3),\theta,x_3)>S_-(\xi(\theta, x_3),\theta,x_3),\ \ \forall (\theta,x_3) \in E.
\ee

The existence and uniqueness of the super-Alfv\'{e}nic and supersonic flow to \eqref{mhd-cyl} follows from the theory of classical solutions to the boundary value problem for quasi-linear symmetric hyperbolic equations (see \cite{BS2007}).
\begin{lemma}\label{supersonic}
{\it Suppose that $\bar{\kappa}^2$ satisfies \eqref{supA1}. Given the incoming data \eqref{super1} satisfying the compatibility conditions
\be\lab{super3}\q\q\begin{cases}
(U_{20},\p_{\theta}^2U_{20},\p_{\theta}(U_{10},U_{30}, P_{0},S_{0},\kappa_{0}))(\pm\theta_0, x_3)=0,\ &\forall x_3\in [-1,1],\\
(U_{30},\p_{x_3}^2 U_{30},\p_{x_3}(U_{10},U_{20}, P_{0},S_{0},\kappa_{0}))(\theta, \pm 1)=0,\  &\forall \theta\in [-\theta_0,\theta_0],
\end{cases}\ee
then there exists $\epsilon_0>0$ depending only on the background solution and the boundary data, such that for any $0<\epsilon<\epsilon_0$, there exists a unique $C^{2,\alpha}(\overline{\Omega})$ solution $(U_{1-}, U_{2-}, U_{3-}, P_-, S_-, \kappa_-)$ to \eqref{mhd-cyl} with \eqref{super1}-\eqref{slip1}, satisfying
\be\no
\|(U_{1-}, U_{2-}, U_{3-}, P_-, S_-, \kappa_-)-(\bar{U}_-,0,0,\bar{P}_-, \bar{S}_-,\bar{\kappa}_-)\|_{C^{2,\alpha}(\overline{\Omega})}\leq C_0\epsilon,
\ee
and
\be\label{super5}\begin{cases}
(U_{2-},\p_{\theta}^2U_{2-},\p_{\theta}(U_{1-},U_{3-}, P_-,S_-,\kappa_-))(r,\pm\theta_0, x_3)=0,&\hspace{-0.5em}\text{on } [r_1,r_2]\times [-1,1],\\
(U_{3-},\p_{x_3}^2U_{3-},\p_{x_3}(U_{1-},U_{2-}, P_-,S_-,\kappa_-))(r,\theta, \pm 1)=0, &\hspace{-0.5em}\text{on } [r_1,r_2]\times [-\theta_0,\theta_0].
\end{cases}\ee
}\end{lemma}

Therefore, our problem is reduced to solving a free boundary value problem for the steady MHD equations in which the shock front and the downstream super-Alfv\'{e}nic subsonic flows are unknown. Then the main result is stated as follows.
\begin{theorem}\label{existence}
{\it Assume that the compatibility conditions \eqref{pressure-cp} and  \eqref{super3} hold and $\bar{\kappa}^2$ satisfies \eqref{supA0}. There exists a suitable constant $\epsilon_0>0$ depending only on the background solutions and the boundary data $U_{10},U_{20},U_{30}, P_{0}, S_{0},\kappa_{0}$, $T_{e}$ such that if $0< \epsilon<\epsilon_0$, the problem \eqref{mhd-cyl} with \eqref{super1}-\eqref{pressure}, and \eqref{rh} has a unique solution $(U_{1+},U_{2+},U_{3+},P_+,S_+,\kappa_+)$ with the shock front $\mlS: r=\xi(\theta,x_3)$ satisfying the following properties.
\begin{enumerate}[(1)]
\item The function $\xi(\theta,x_3)\in C^{3,\alpha}(\overline{E})$ satisfies
\be\no
\|\xi(\theta, x_3)-r_s\|_{C^{3,\alpha}(\overline{E})}\leq C_*\epsilon,
\ee
and
\be\no\begin{cases}
\p_{\theta}\xi(\pm\theta_0, x_3)=\p_{\theta}^3\xi(\pm\theta_0, x_3)=0,\ \ &\forall x_3\in [-1,1],\\
\p_{x_3}\xi(\theta,\pm 1)=\p_{x_3}^3\xi(\theta, \pm 1)=0,\ \ &\forall \theta\in [-\theta_0,\theta_0],
\end{cases}\ee
where $C_*$ is a positive constant depending only on the background solution, the supersonic incoming flow, and the exit pressure.
\item The solution $(U_{1+},U_{2+},U_{3+},P_+,S_+,\kappa_+)\in C^{2,\alpha}(\overline{\Omega_+})$ satisfies the entropy condition
\begin{equation}\no
S_+(\xi(\theta,x_3)+,\theta,x_3) >S_-(\xi(\theta,x_3)-, \theta,x_3)\quad \forall (\theta,x_3)\in E,
\end{equation}
and the estimate
\be\no
\|(U_{1+}, U_{2+}, U_{3+}, P_+, S_+, \kappa_+)-(\bar{U}_+,0,0,\bar{P}_+, \bar{S}_+,\bar{\kappa}_+)\|_{C^{2,\alpha}(\overline{\Omega_+})}\leq C_*\epsilon,
\ee
with the compatibility conditions
\be\no
\begin{cases}
(U_{2+},\p_{\theta}^2U_{2+},\p_{\theta}(U_{1+},U_{3+}, P_+,S_+,\kappa_+))(r,\pm\theta_0, x_3)=0,\hspace{-0.5em}&\forall (r,x_3)\in [\xi,r_2]\times [-1,1],\\
(U_{3+},\p_{x_3}^2U_{3+},\p_{x_3}(U_{1+},U_{2+}, P_+,S_+,\kappa_+))(r,\theta, \pm 1)=0,\hspace{-0.5em}&\forall (r,\theta)\in [\xi,r_2]\times [-\theta_0,\theta_0].
\end{cases}\ee
\end{enumerate}
}\end{theorem}

This paper will be organized as follows. In Section \ref{reformulation}, we decompose the hyperbolic and elliptic modes for the steady MHD equations in terms of the deformation and the modified vorticity, and reformulate the Rankine-Hugoniot jump conditions, which are well-suited to our decomposition of the steady MHD equations. In Section \ref{proof}, we design an iteration scheme and solve the deformation-curl system with nonlocal terms and the unusual second order differential boundary condition on the shock front. 

\section{The reformulation of the shock problem}\label{reformulation}

\subsection{The reformulation of the steady MHD equations \texorpdfstring{\eqref{mhd-cyl}}{} in terms of the deformation tensor and the modified vorticity}

In the downstream subsonic region, we will show that the steady MHD equations \eqref{mhd-cyl} are elliptic-hyperbolic mixed if the magnetic field satisfies the assumption \eqref{supA2}. First, it is easy to identify the hyperbolic modes in \eqref{mhd-cyl}. The Bernoulli's quantity, the entropy, and the scalar function $\kappa$ are conserved along the streamlines:
\be\label{ber10}
(\p_r+\frac{U_2}{U_1}\frac{1}{r}\p_{\theta} +\frac{U_3}{U_1}\p_{x_3}) (B,S,\kappa)=0.
\ee

To go further, we move back to the steady MHD equations \eqref{mhd1} in Cartesian coordinates. Using the vector identity ${\bf u}\cdot\nabla {\bf u}=\text{curl }{\bf u}\times {\bf u}+\nabla \frac12 |{\bf u}|^2$, the momentum equations can be rewritten as
\be\no
&&\text{curl }{\bf u}\times {\bf u} + \nabla \frac12 |{\bf u}|^2 + \frac{1}{\rho}\nabla P =\frac{1}{\rho}\text{curl }{\bf h}\times {\bf h}=\kappa \text{curl }(\kappa \rho {\bf u})\times {\bf u}\\\no
&&=\text{curl }(\kappa^2 \rho {\bf u})\times {\bf u}-(\nabla\kappa\times (\kappa \rho {\bf u})\times {\bf u}.
\ee
Therefore
\be\label{mo2}
\text{curl }\{(1-\kappa^2 \rho){\bf u}\}\times {\bf u}+ \nabla B-\frac{\rho^{\gamma-1}}{\gamma-1} \nabla S=-(\nabla\kappa\times (\kappa \rho {\bf u})\times {\bf u}=\kappa \rho |{\bf u}|^2 \nabla \kappa.
\ee

Motivated by \eqref{mo2}, we introduce the modified vorticity ${\bf J}=\text{curl }\{(1-\kappa^2\rho){\bf u}\}=J_1 {\bf e}_r + J_2 {\bf e}_{\theta}+ J_3 {\bf e}_3$, where
\be\no\begin{cases}
J_1= \frac{1}{r}\p_{\theta} \{(1-\kappa^2\rho)U_3\}- \p_{x_3} \{(1-\kappa^2\rho)U_2\},\\
J_2=\p_{x_3} \{(1-\kappa^2\rho)U_1\}- \p_r \{(1-\kappa^2\rho)U_3\},\\
J_3=(\p_r+\frac1r) \{(1-\kappa^2\rho)U_2\}- \frac{1}{r}\p_\theta \{(1-\kappa^2\rho)U_1\}.
\end{cases}\ee
Clearly, the modified vorticity records the effect of the Lorentz force on the velocity field.

Transforming the equations \eqref{mo2} to the cylindrical coordinates, one obtains that
\be\no\begin{cases}
U_1 J_3 - U_3 J_1 + \frac{1}{r}\p_{\theta} B-\frac{B-\frac{1}{2}|{\bf U}|^2} {\gamma S}\frac{1}{r}\p_{\theta}S-\kappa \rho |{\bf U}|^2 \frac1r\p_{\theta}\kappa=0,\\
-U_1 J_2 + U_2 J_1 +\p_{x_3} B-\frac{B-\frac{1}{2}|{\bf U}|^2} {\gamma S U_1}\p_{x_3} S+\kappa \rho |{\bf U}|^2\p_{x_3}\kappa=0.
\end{cases}\ee
Thus there holds
\be\label{vor13}\begin{cases}
J_2 =\frac{1}{U_1}(U_2 J_1+\p_{x_3} B- \frac{B-\frac{1}{2}|{\bf U}|^2} {\gamma S}\p_{x_3}S-\kappa \rho |{\bf U}|^2 \p_3 \kappa),\\
J_3 =\frac1{U_1}(U_3 J_1-\frac{1}{r}\p_{\theta} B + \frac{B-\frac{1}{2}|{\bf U}|^2} {\gamma S}\frac{1}{r}\p_{\theta} S+\kappa \rho |{\bf U}|^2 \frac1{r}\p_{\theta}\kappa).
\end{cases}\ee

Since
\be\no
\text{div }\text{curl }{\{(1-\kappa^2\rho)\bf u\}}=\p_r J_1+\frac{1}{r}\p_{\theta} J_2+\p_{x_3} J_3 +\frac{1}{r} J_1=0,
\ee
substituting \eqref{vor13} into the above equation yields
\be\label{vor14}
&&(\p_r+\frac{U_2}{U_1}\f{1}{r}\p_{\theta} +\frac{U_3}{U_1}\p_{x_3})J_1 +(\f{1}{r} +\frac{1}{r}\p_{\theta}(\f{U_2}{U_1})+\p_{x_3}(\frac{U_3}{U_1}))J_1\\\no
&&\quad+ \frac{1}{r}\p_{\theta}(\frac1{U_1}) \p_{x_3}B-\p_{x_3}(\f{1}{U_1}) \frac{1}{r}\p_{\theta} B-\f{1}{r} \p_{\theta}(\frac{B-\frac{1}{2}|{\bf U}|^2}{\gamma S U_1}) \p_{x_3} S \\\no&&
\quad+ \p_{x_3} (\frac{B-\frac{1}{2}|{\bf U}|^2}{\gamma S U_1}) \frac{1}{r}\p_{\theta} S-\frac{1}{r}\p_{\theta}(\frac{\kappa \rho |{\bf U}|^2}{U_1})\p_{x_3} \kappa+\p_{x_3}(\frac{\kappa \rho |{\bf U}|^2}{U_1})\frac1{r}\p_{\theta} \kappa=0.
\ee

Next, we study the elliptic modes in the steady MHD equations \eqref{mhd-cyl}. Using the Bernoulli's quantity $B=\frac12 |{\bf U}|^2 + \frac{\gamma P}{(\gamma-1)\rho}$, one can represent the density as a function of $B, S$, and $|{\bf U}|^2$:
\be\lab{den}
\rho= \rho(B,S,|{\bf U}|^2)=\b(\frac{\gamma-1}{\gamma S}\b)^{\frac{1}{\gamma-1}}\b(B-\frac12 |{\bf U}|^2\b)^{\frac{1}{\gamma-1}}.
\ee
Substituting \eqref{den} into the continuity equation and using \eqref{ber10} lead to
\be\no
&&(c^2(B,|{\bf U}|^2)-U_1^2)\p_r U_1 + (c^2(B,|{\bf U}|^2)-U_2^2)\frac{1}{r}\p_{\theta} U_2+ (c^2(B,|{\bf U}|^2)-U_3^2)\p_{x_3} U_3\\\no
&&\quad
+ \frac{c^2(B,|{\bf U}|^2) U_1}{r}=U_1(U_2\p_r U_2+ U_3 \p_r U_3)+ U_2(U_1\frac{1}{r}\p_{\theta} U_1+U_3\frac{1}{r}\p_{\theta} U_3)\\\lab{den11}
&&\quad\quad+ U_3 (U_1\p_{x_3} U_1+U_2 \p_{x_3} U_2),
\ee
which can be rewritten as a Frobenius inner product of a symmetric matrix and the deformation matrix.

The equation \eqref{den11} together with the vorticity equations constitutes a deformation-curl system for the velocity field:
\be\label{dc}\begin{cases}
(c^2-U_1^2)\p_r U_1 + (c^2-U_2^2)\frac{1}{r}\p_{\theta} U_2\\
\q + (c^2-U_3^2)\p_{x_3} U_3 + \frac{c^2 U_1}{r}=U_1(U_2\p_r U_2+ U_3 \p_r U_3)\\
\q + U_2(U_1\frac{1}{r}\p_{\theta} U_1+U_3\frac{1}{r}\p_{\theta} U_3)+ U_3 (U_1\p_{x_3} U_1+U_2 \p_{x_3} U_2),\\
\frac{1}{r}\p_{\theta} \{(1-\kappa^2\rho) U_3\}- \p_{x_3}\{(1-\kappa^2\rho) U_2\}=J_1,\\
\p_{x_3} \{(1-\kappa^2\rho)U_1)- \p_r \{(1-\kappa^2\rho)U_3\}=J_2,\\
(\p_r+\frac1r) \{(1-\kappa^2\rho)U_2)- \frac{1}{r}\p_\theta \{(1-\kappa^2\rho)U_1\}=J_3.
\end{cases}\ee

\begin{lemma}\label{equiv}({\bf Equivalence.})
{\it Assume that $C^2$ smooth vector functions $(\rho, {\bf U}, S, \kappa)$ defined on $\Omega$ do not contain the vacuum (i.e., $\rho>0$ in $\Omega$) and the radial velocity $U_1$ is always positive in $\Omega$, then the following two statements are equivalent:
\begin{enumerate}[(1).]
\item $(\rho, {\bf U}, S,\kappa)$ satisfy the steady MHD system \eqref{mhd-cyl} in $\Omega$;
\item $({\bf U}, S, B,\kappa)$ satisfy the equations \eqref{ber10}, \eqref{vor13} and \eqref{dc}.
\end{enumerate}
}\end{lemma}

To simplify the notation, we set
\be\no
&&W_{1}(r,\theta,x_3)= U_{1+}(r,\th,x_3)- \bar{U}_+(r),\quad W_{j}(r,\theta,x_3)=U_{j+}(r,\theta,x_3), j=2,3,\\\no
&&W_{4}(r,\theta,x_3)=S_+(r,\theta,x_3)-\bar{S}_+,\ W_{5}(r,\theta,x_3)= B_+(r,\theta,x_3)-\bar{B},\ \\\no
&&W_{6}(r,\theta,x_3)=\kappa_{+}(r,\theta,x_3)-\bar{\kappa},\ W_7(\theta,x_3)= \xi(\theta,x_3)-r_s, \ \ \ {\bf W}=(W_1,\cdots, W_6).
\ee
Then the density and the pressure can be expressed as
\be\label{denw1}
&&\rho({\bf W})=\b(\frac{\gamma-1}{\gamma (\bar{S}+W_4)}\b)^{\frac{1}{\gamma-1}}\b(\bar{B}+W_5-\frac12(\bar{U}+W_1)^2-\frac12\sum_{i=2}^3 W_i^2\b)^{\frac{1}{\gamma-1}},\\\label{denw2}
&&P({\bf W})=\b(\frac{(\gamma-1)^{\gamma}}{\gamma^{\gamma}(\bar{S}+W_4)}\b)^{\frac{1}{\gamma-1}}\b(\bar{B}+W_5-\frac12(\bar{U}+W_1)^2-\frac12\sum_{i=2}^3 W_i^2\b)^{\frac{\gamma}{\gamma-1}}.
\ee
Before we write the equations \eqref{ber10}, \eqref{vor13}, \eqref{vor14} and \eqref{dc} in terms of ${\bf W}$, we introduce some notations.
\be\no
&&d_1(r)=1-\bar{M}_+^2(r), \ d_2(r)=\frac{\bar{M}_+^2(2+(\gamma-1)\bar{M}_+^2)}{r(1-\bar{M}_+^2)},\\\no
&& d(r)=1-\bar{\ka}^2\bar{\rho}_+(r)+\bar{\ka}^2\bar{\rho}_+(r)\bar{M}_+^2(r), \ \ d_0(r)=1-\bar{\ka}^2\bar{\rho}_+(r),\\\no
&&d_{3}(r)=\frac{\bar{B}-\frac12 \bar{U}_+^2(r)}{\gamma \bar{S}_+ \bar{U}_+(r)}+\frac{\bar{\ka}^2(\bar{\rho}\bar{U})_+(r)}{(\gamma-1)\bar{S}_+},\\\no
&& \bar{c}_+^2(r)=c^2(\bar{B},\bar{U}_+^2(r)), \ \ d_5(r)=-\frac{(\gamma-1)(\bar{U}_+'+\frac{\bar{U}_+}{r})}{\bar{c}_+^2(r)}.
\ee

The equations for the hyperbolic quantities $W_4, W_5$ and $W_6$ are
\be\label{ent11}
\b(\p_r+\frac{W_2}{\bar{U}+W_1}\frac{1}{r}\p_{\theta}+\frac{W_3}{\bar{U}+W_1}\p_{x_3}\b) (W_4,W_5,W_6)=0,
\ee
In terms of ${\bf W}$, the vorticity ${\bf J}$ has the form
\be\no\begin{cases}
J_1=\frac1r\p_{\theta}\{(1-(\bar{\ka}+W_6)^2\rho({\bf W}))W_3\}-\p_{x_3}\{(1-(\bar{\ka}+W_6)^2\rho({\bf W}))W_2\},\\
J_2=\p_{x_3}\{(1-(\bar{\ka}+W_6)^2\rho({\bf W}))W_1\}-\p_{r}\{(1-(\bar{\ka}+W_6)^2\rho({\bf W}))W_3\}\\
\q\q-\bar{U}\p_{x_3}\{(\bar{\ka}+W_6)^2\rho({\bf W}))\},\\
J_3 =(\p_r+\frac1r)\{(1-(\bar{\ka}+W_6)^2\rho({\bf W}))W_2\}-\frac1r\p_{\theta}\{(1-(\bar{\ka}+W_6)^2\rho({\bf W}))W_1\}\\
\q\q\q\q+\frac{\bar{U}}{r}\p_{\theta}\{(\bar{\ka}+W_6)^2\rho({\bf W})\}.
\end{cases}\ee
The equations for the vorticity $J_1$ are
\be\no
&&\b(\p_r+\frac{W_2}{\bar{U}+W_1}\frac{1}{r}\p_{\theta} +\frac{W_3}{\bar{U}+W_1}\p_{x_3}\b)J_1 +\b(\f{1}{r} +\frac{1}{r}\p_{\theta}\b(\frac{W_2}{\bar{U}+W_1}\b)+\p_{x_3}\b(\frac{W_3}{\bar{U}+W_1}\b)\b)J_1\\\no
&&\quad\quad+ \frac{1}{r}\p_{\theta}\b(\frac1{\bar{U}+W_1}\b) \p_{x_3}W_5-\p_{x_3}\b(\frac{1}{\bar{U}+W_1}\b) \frac{1}{r}\p_{\theta} W_5\\\no
&&\quad\quad-\f{1}{r} \p_{\theta} \b(\frac{\bar{B}-\frac{1}{2}\bar{U}^2+W_5-\bar{U}W_1-\frac{1}{2}\sum_{i=1}^3W_i^2}{\gamma (\bar{S}+W_4)(\bar{U}+W_1)}\b) \p_{x_3} W_4\\\label{vor21}
&&\quad\quad+ \p_{x_3} \b(\frac{\bar{B}-\frac{1}{2}\bar{U}^2+W_5-\bar{U}W_1-\frac{1}{2}\sum_{i=1}^3W_i^2}{\gamma (\bar{S}+W_4)(\bar{U}+W_1)}\b) \frac{1}{r}\p_{\theta} W_4\\\no
&&\quad\quad-\frac1{r}\p_{\theta}\b(\frac{(\bar{\kappa}+W_6)\rho({\bf W})((\bar{U}+W_1)^2+W_2^2+W_3^2)}{\bar{U}+W_1}\b)\p_{x_3}W_6\\\no
&&\quad\quad+\p_{x_3}\b(\frac{(\bar{\kappa}+W_6)\rho({\bf W})((\bar{U}+W_1)^2+W_2^2+W_3^2)}{\bar{U}+W_1}\b)\frac1r\p_{\theta}W_6=0,
\ee
and
\be\label{vor22}\begin{cases}
J_2=\frac{W_2J_1 +\p_{x_3} W_5-(\bar{\kappa}+W_6)\rho({\bf W})((\bar{U}+W_1)^2+W_2^2+W_3^2)\p_{x_3}W_6}{\bar{U}+W_1}\\
\q\q\q- \frac{\bar{B}-\frac{1}{2}\bar{U}^2+W_5-\bar{U}W_1-\frac{1}{2}\sum_{i=1}^3W_i^2} {\gamma (\bar{S}+W_4)(\bar{U}+W_1)}\p_{x_3} W_4,\\
J_3=\frac{W_3J_1 -\frac{1}{r}\p_{\theta} W_5+ (\bar{\kappa}+W_6)\rho({\bf W})((\bar{U}+W_1)^2+W_2^2+W_3^2)\frac1r\p_{\theta} W_6}{\bar{U}+W_1} \\
\q\q\q+ \frac{\bar{B}-\frac{1}{2}\bar{U}^2+W_5-\bar{U}W_1-\frac{1}{2}\sum_{i=1}^3W_i^2} {\gamma (\bar{S}+W_4)(\bar{U}+W_1)}\frac{1}{r}\p_{\theta} W_4.
\end{cases}\ee

It follows from \eqref{den11} that
\be\label{den15}
d_1(r)\p_r W_1 + \frac{1}{r}\p_{\theta} W_2 +\p_{x_3} W_3 + (\frac1{r}+d_2(r))W_1= d_5(r) W_5+ F({\bf W}),
\ee
where
\be\no
&&\bar{c}^2(r) F({\bf W})=-(\gamma-1)(\p_r W_1+\frac{W_1}{r}) W_5\\\no
&&\quad\quad+ (\bar{U}'+\p_r W_1)(\frac{\gamma+1}{2} W_1^2+\frac{\gamma-1}{2}(W_2^2+W_3^2))\\\no
&&\quad\quad+\frac{(\gamma-1)(\bar{U}+W_1)}{2r }\sum_{j=1}^3 W_j^2+(\gamma+1)\bar{U}W_1 \p_r W_1+ (\gamma-1)\bar{U} \frac{W_1^2}{r}\\\no
&&\quad\quad-\{(\gamma-1)W_5-\frac{\gamma-1}{2}\sum_{j=1}^3 W_j^2-(\gamma-1)\bar{U}W_1\}(\frac{1}{r}\p_{\theta} W_2+\p_{x_3} W_3)\\\no
&&\quad\quad
+(W_2^2\frac{1}{r}\p_{\theta} W_2 +W_3^2 \p_{x_3}W_3)+ (\bar{U}+W_1)(W_2\p_r W_2+W_3\p_r W_3)\\\no
&&\quad\quad+W_2((\bar{U}+W_1)\frac{1}{r}\p_{\theta} W_1 +\frac{W_3}{r}\p_{\theta} W_3)+ W_3((\bar{U}+W_1)\p_{x_3} W_1 +W_2\p_{x_3} W_2).
\ee
Here $F({\bf W})$ and the following $H_i, G_i$ are quadratic and high order terms.

\subsection{The linearization of the Rankine-Hugoniot conditions and boundary conditions}\label{22}

For the super-Alfv\'{e}nic and subsonic flow, the steady MHD system is elliptic-hyperbolic mixed, and thus it is important to formulate proper boundary conditions and their compatibility.

It follows from the third and fourth equations in \eqref{rh} that
\be\lab{shock11}
\frac{1}{\xi(\theta,x_3)}\p_{\theta} \xi =\frac{f_2(\xi,\theta,x_3)}{f(\xi,\theta,x_3)},\ \,\,\, \p_{x_3}  \xi =\frac{f_3(\xi,\theta,x_3)}{f(\xi,\theta,x_3)},
\ee
where
\be\no\begin{cases}
f(\xi,\theta,x_3)=   [\rho U_2^2 + P+\frac12 \kappa^2 \rho^2 (U_1^2+U_3^2-U_2^2)] [\rho U_3^2 +P+\frac12 \kappa^2 \rho^2 (U_1^2+U_2^2-U_3^2)]\\
\q\q- ([(1-\kappa^2\rho)\rho U_2 U_3])^2,\\
f_2(\xi,\theta,x_3)= [\rho U_3^2+P+\frac12 \kappa^2 \rho^2 (U_1^2+U_2^2-U_3^2)] [(1-\kappa^2\rho)\rho U_1 U_2]\\
\q\q- [(1-\kappa^2\rho)\rho U_1 U_3] [(1-\kappa^2\rho)\rho U_2 U_3],\\
f_3(\xi,\theta,x_3)= [\rho U_2^2+P+\frac12 \kappa^2 \rho^2 (U_1^2+U_3^2-U_2^2)][(1-\kappa^2\rho)\rho U_1 U_3]\\
\q\q- [(1-\kappa^2\rho)\rho U_1 U_2] [(1-\kappa^2\rho)\rho U_2 U_3].
\end{cases}\ee
Rewrite \eqref{shock11} as
\be\lab{shock12}\q\q\begin{cases}
\p_{\theta} \xi= a_0 r_sW_2 +r_s g_2(\bm{\Psi}_-(r_s+W_7,\theta,x_3) - \overline{\bm{\Psi}}_-(r_s+W_7),{\bf W}, W_7),\\
\p_{x_3} \xi= a_0 W_3 +g_3(\bm{\Psi}_-(r_s+W_7,\theta,x_3) - \overline{\bm{\Psi}}_-(r_s+W_7),{\bf W}, W_7),
\end{cases}\ee
where $a_0= d_0(r_s)\f{\bar{\rho}_+ \bar{U}_+}{[\bar{P}]}(r_s)>0$ and
\be\no
&&g_2=\frac{1}{r_s}\b(\frac{\xi f_2}{f}-r_s a_0  d_0(r_s)W_2(\xi(\th,x_3),\th,x_3)\b),\ \ \\\no
&&g_3=\frac{f_3}{f}- a_0 d_0(r_s)W_3(\xi(\th,x_3),\th,x_3).
\ee
The functions $g_i, i=2,3$ are error terms which can be bounded by
\be\label{g23}
|g_i|\leq C_*(|\bm{\Psi}_-(r_s+W_7,\theta,x_3) - \overline{\bm{\Psi}}_-(r_s+W_7)|+|{\bf W}|^2+|W_7|^2).
\ee

It follows from \eqref{shock11} and \eqref{rh} that
\be\label{rh1}
\begin{cases}
[\rho U_1]=\frac{[\rho U_2] f_2+[\rho U_3] f_3}{f},\\
[\rho U_1^2+P]=\frac{1}{f}\sum_{i=2}^3\{[(1-\kappa^2\rho)\rho U_1 U_i]\\
+\frac12 \kappa^2 (\rho_+ U_{1+} + \rho_- U_{1-})[\rho U_i]\}f_i-\frac12 [\kappa^2 \rho^2(U_2^2+U_3^2)],\\
[B]=[\kappa]=0.
\end{cases}\ee

Note that
\be\no
[\bar{\rho} \bar{U}](r_s+W_7)=O(W_7^2),\q [\bar{\rho}\bar{U}^2+\bar{P}](r_s+W_7)=-\frac{1}{r_s}[\bar{P}(r_s)] W_7+O(W_7^2).
\ee

By the Taylor's expansion and \eqref{rh1}, it holds that at $(\xi(\theta,x_3),\theta,x_3)$
\be\label{rh2}\begin{cases}
a_{11} W_{1}+ a_{12} W_{4} = R_{01}(\bm{\Psi}_-(r_s+W_7,\theta,x_3) - \overline{\bm{\Psi}}_-(r_s+W_7),{\bf W}, W_7),\\
a_{21} W_{1}+ a_{22} W_{4} =-\frac{1}{r_s}[\bar{P}(r_s)] W_7+R_{02}(\bm{\Psi}_-(r_s+W_7,\theta,x_3) - \overline{\bm{\Psi}}_-(r_s+W_7),{\bf W}, W_7),
\end{cases}\ee
where
\be\no
&&a_{11}=(\bar{\rho}_+ (1-\bar{M}_+^2))(r_s),\ \ a_{12}=-\frac{(\bar{\rho}_+\bar{U}_+)(r_s)}{(\gamma-1)\bar{S}_+},\\\no
&&a_{21}= (\bar{\rho}_+\bar{U}_+(1-\bar{M}_+^2))(r_s),\ \ a_{22}=-\b(\frac{(\bar{\rho}_+\bar{U}_+^2)(r_s)}{(\gamma-1)\bar{S}_+}+\frac{1}{\gamma-1}\bar{\rho}_+^{\gamma}(r_s)\bigg),\\\no
&&R_{01}=\frac{[\rho U_2] f_2+[\rho U_3] f_3}{f}-[\bar{\rho}\bar{U}(\xi)]+(\rho_- U_{1-})(\xi,\theta,x_3)-(\bar{\rho}_-\bar{U}_-)(\xi)\\\no
&&\quad-\rho({\bf W})(\bar{U}_+ + W_{1})+(\bar{\rho}_+\bar{U}_+)(\xi)+ a_{11} W_{1}+ a_{12} W_{4},\\\no
&&R_{02}=\frac{1}{f}\sum_{i=2}^3\{[(1-\kappa^2\rho)\rho U_1 U_i] +\frac12 \kappa^2 (\rho_+ U_{1+} + \rho_- U_{1-})[\rho U_i]\}f_i\\\no
&&\q-\frac12 [\kappa^2 \rho^2(U_2^2+U_3^2)]-[(\bar{\rho}\bar{U}^2+\bar{P})(\xi)]+(\rho_- U_{1-}^2+ P_-)(\xi,\theta,x_3)\\\no
&&\quad-(\bar{\rho}_-\bar{U}_-^2+ \bar{P}_-)(\xi)-\rho({\bf W})(\bar{U}_+ + W_{1})^2+ P({\bf W})\\\no
&&\quad+(\bar{\rho}_+\bar{U}_+^2+ \bar{P}_+)(\xi)+ a_{21} W_{1}+ a_{22} W_{4}.
\ee

Then solving the algebraic equations in \eqref{rh2}, one gets at $(\xi,\theta,x_3)$
\be\lab{shock13}\begin{cases}
W_1=a_1 W_7+ R_1(\bm{\Psi}_-(r_s+W_7,\theta,x_3) - \overline{\bm{\Psi}}_-(r_s+W_7),{\bf W}, W_7),\\
W_4=a_2 W_7+R_2(\bm{\Psi}_-(r_s+W_7,\theta,x_3) - \overline{\bm{\Psi}}_-(r_s+W_7),{\bf W}, W_7),\\
W_5=B_-(r_s+W_7(\theta,x_3),\theta,x_3)-\bar{B}_-,\\
W_6=\kappa_-(r_s+W_7(\theta,x_3),\theta,x_3)-\bar{\kappa},
\end{cases}\ee
where
\be\no
&&a_1=\f{\gamma \bar{U}_+(r_s)[\bar{P}(r_s)]}{r_s\bar{\rho}_+(r_s) (c^2(\bar{\rho}_+(r_s),\bar{S}_+)-\bar{U}_+^2(r_s))}>0,\\\no
&&a_2=\frac{(\gamma-1)[\bar{P}(r_s)]}{r_s\bar{\rho}_+^{\gamma}(r_s)}>0,
\ee
and
\begin{equation*}
\begin{split}
R_1&=\frac{a_{12}R_{02}-a_{22}R_{01}}{a_{11}a_{22}-a_{12}a_{21}}:=\sum_{i=1}^{2}b_{1i}R_{0i},\\
R_2&=\frac{a_{21}R_{01}-a_{11}R_{02}}{a_{11}a_{22}-a_{12}a_{21}}:=\sum_{i=1}^{2}b_{2i}R_{0i}.
\end{split}
\end{equation*}
There exists $C_0>0$ depending only on the background solution, such that
\be\no
|R_{i}|\leq C_0(|\bm{\Psi}_-(r_s+W_7,\theta,x_3) - \overline{\bm{\Psi}}_-(r_s+W_7)|+ |{\bf W}(\xi,\theta,x_3)|^2+W_7^2), i=1,2.
\ee

It follows from \eqref{denw1}-\eqref{denw2} and the Taylor's expansion that
\be\no
&&\epsilon T_e(\theta,x_3)=-\bar{\rho}_+\bar{U}_+(r_2)(d(r_2)W_1 + d_{3}(r_2)W_4)(r_2,\theta,x_3)\\\no
&&\q\q\q+(\bar{\rho}+\bar{\ka}^2 \bar{\rho}_+^2 \bar{M}_+^2) W_5+ \bar{\ka}\bar{\rho}^2\bar{U}^2 W_6+E({\bf W}(r_2,\theta,x_3)),
\ee
where
\be\no
&&E({\bf W})(r_2,\cdot)=P({\bf W})+ \frac12 (\bar{\kappa}+W_{6})^2(\rho({\bf W}))^2((\bar{U}_++W_{1})^2 + W_{2}^2 +W_{3}^2)(r_2,\cdot)\\\no
&&\quad-(\bar{P}_++\frac12 \bar{\kappa}^2 \bar{\rho}_+^2 \bar{U}_+^2)(r_2)+ (\bar{\rho}_+\bar{U}_+)(r_2)(d(r_2) W_1 + d_{3}(r_2) W_{4})(r_2,\cdot)\\\no
&&\quad -\bar{\rho}_+(r_2)(1+\bar{\kappa}^2 \bar{\rho}_+ \bar{M}_+^2)(r_2)W_5(r_2,\cdot)-\bar{\kappa}(\bar{\rho}_+\bar{U}_+)(r_2) W_6(r_2,\cdot).
\ee
and $E$ is an error term that can be bounded by
\be\no
|E({\bf W}(r_2,\theta,x_3))|\leq C_*|{\bf W}(r_2,\theta,x_3)|^2.
\ee

This, together with \eqref{pressure} implies that at $(r_2,\theta,x_3)$ there holds
\be\no
&&\{d W_{1} + d_{3} W_{4}\}(r_2,\cdot)=-\frac{\epsilon T_e(\cdot)}{(\bar{\rho}_+\bar{U}_+)(r_2)}+\frac{(1+\bar{\kappa}^2 \bar{\rho}_+ \bar{M}_+^2(r_2))W_{5}(r_2,\cdot)}{\bar{U}_+(r_2)}\\\label{pres2}
&&\quad\quad-\bar{\kappa}(\bar{\rho}_+\bar{U}_+)(r_2) W_6(r_2,\cdot)-\frac{E({\bf W}(r_2,\theta,x_3))}{(\bar{\rho}_+\bar{U}_+)(r_2)}.
\ee

The boundary conditions for $W_{2}$ and $W_{3}$ on the nozzle walls are
\be\label{slip15}\begin{cases}
W_{2}(r,\pm\theta_0, x_3)=0,\ \ &\text{on } r_s+W_7(\pm\theta_0,x_3)<r<r_2,x_3\in [-1,1],\\
W_{3}(r,\theta,\pm 1)=0,\ \ &\text{on }r_s+W_7(\theta,\pm 1)<r<r_2,\theta\in [-\theta_0,\theta_0].
\end{cases}\ee

Then to solve the problem \eqref{mhd-cyl} with \eqref{super1}-\eqref{pressure}, and \eqref{rh} is equivalent to finding a function $W_7$ defined on $E$ and vector functions $(W_1,\cdots, W_6)$ defined on the $\Omega_{W_7}:=\{(r,\theta,x_3): r_s +W_7(\theta,x_3)<r<r_2,(\theta,x_3)\in E\}$, which solves the equations \eqref{ent11}--\eqref{den15} with boundary conditions \eqref{shock12},\eqref{shock13},\eqref{pres2} and \eqref{slip15}.

In the following, the subscript ``+" will be ignored to simplify the notations.

\subsection{Transform to a fixed boundary value problem} To fix the subsonic region, relabeling $V_7(\theta,x_3)=\xi(\theta,x_3)-r_s$, we introduce the coordinates transformation
\be\label{coor}
y_1=\frac{r-r_s-V_7}{r_2-r_s-V_7}(r_2-r_s) + r_s,\ y_2=\theta,\  y_3=x_3.
\ee
Then
\be\no\begin{cases}
r= y_1+\frac{r_2-y_1}{r_2-r_s}V_7=: D_0^{V_7},\\
\p_r=\frac{r_2-r_s}{r_2-r_s-V_7(y_2,y_3)} \p_{y_1}=: D_1^{V_7},\\
\frac{1}{r}\p_{\theta}=\frac{1}{D_0^{V_7}}(\p_{y_2}+\frac{(y_1-r_2)\p_{y_2}V_7}{r_2-r_s-V_7}\p_{y_1})=: D_2^{V_7},\\
\p_{x_3}=\p_{y_3}+\frac{(y_1-r_2)\p_{y_3}V_7}{r_2-r_s-V_7}\p_{y_1}=: D_3^{V_7},
\end{cases}\ee
and the domain $\Omega_+$ is changed to be
\be\no
\mathbb{D}=\{(y_1,y'): y_1\in (r_s, r_2), y'=(y_2,y_3)\in E\}.
\ee
Denote
\be\no
&&\Sigma_2^{\pm}=\{(y_1,\pm\theta_0,y_3):(y_1,y_3)\in (r_s,r_2)\times (-1,1)\},\\\no
&&\Sigma_3^{\pm}=\{(y_1,y_2,\pm 1):(y_1,y_2)\in (r_s,r_2)\times (-\theta_0,\theta_0)\}.
\ee

Set
\be\no\begin{cases}
V_i(y)= W_i(y_1+\f{r_2-y_1}{r_2-r_s}V_7,y'), i=1,\cdots, 6,\\
\tilde{J}_i(y)=J_i(y_1+\f{r_2-y_1}{r_2-r_s}V_7,y'), i=1,2,3.
\end{cases}\ee
Then the functions $\rho(r,\theta,x_3)$ and $P(r,\theta,x_3)$ in \eqref{denw1}-\eqref{denw2} are transformed to be
\be\no
&&\tilde{\rho}({\bf V}(y),V_7)=\bigg(\frac{\gamma-1}{\gamma(\bar{S}+V_4)}\bigg)^{\frac{1}{\gamma-1}}\bigg(\bar{B}+V_5-\frac{1}{2}(\bar{U}(D_0^{V_7})+V_1)^2-\frac{1}{2}\sum_{i=2}^3 V_i^2\bigg)^{\frac{1}{\gamma-1}},\\\no
&&\tilde{P}({\bf V}(y),V_7)=\bigg(\frac{(\gamma-1)^{\gamma}}{\gamma^{\gamma}(\bar{S}+V_4)}\bigg)^{\frac{1}{\gamma-1}}\bigg(\bar{B}+V_5-\frac{1}{2}(\bar{U}(D_0^{V_7})+V_1)^2-\frac{1}{2}\sum_{i=2}^3 V_i^2\bigg)^{\frac{\gamma}{\gamma-1}}.
\ee
In the $y$-coordinates, \eqref{shock12} is changed to be
\begin{equation}\label{D23_V7_ycor}
\begin{split}
\frac{1}{r_s}\p_{y_2} V_7(y')&=a_0V_2(r_s,y') + g_2({\bf V}(r_s,y'), V_7),\\
\p_{y_3} V_7(y')&=a_0 V_3(r_s,y') + g_3({\bf V}(r_s,y'), V_7),
\end{split}
\end{equation}
where
\be\label{g21}
&&g_2=\frac{1}{r_s}\b(\frac{(r_s+V_7)f_2({\bf V}(r_s,y'), V_7(y'))}{f({\bf V}(r_s,y'), V_7(y'))}-a_0 r_s V_2(r_s,y')\b),\\\label{g31}
&&g_3=\frac{1}{r_s}\b(\frac{f_3({\bf V}(r_s,y'), V_7(y'))}{f({\bf V}(r_s,y'), V_7(y'))}-a_0 V_3(r_s,y')\b).
\ee

In the $y$ coordinates, the transonic shock problem can be reformulated as follows. The shock front will be determined by the first equation in \eqref{shock13} as follows:
\be\label{shock400}
V_7(y')=a_1^{-1}V_1(r_s,y')- a_1^{-1}R_1({\bf V}(r_s,y'),V_7(y')),
\ee
where $R_{1}({\bf V}(r_s,y'),V_7(y'))=\displaystyle\sum_{i=1}^2 b_{1i}R_{0i}({\bf V}(r_s,y'),V_7(y'))$.

The last three formulas in \eqref{shock13} will be used to solve the hyperbolic modes:
\be\label{ber31}\begin{cases}
\b(D_1^{V_7}+\frac{V_2}{\bar{U}(D_0^{V_7})+V_1}D_2^{V_7}+\frac{V_3}{\bar{U}(D_0^{V_7})+V_1}D_3^{V_7}\b)(V_4,V_5,V_6)=0,\\
V_4(r_s,y')=a_2 V_7(y')+R_2({\bf V}(r_s,y'),V_7(y')),\\
V_5(r_s,y')=B_-(r_s+V_7(y'),y')-\bar{B},\\
V_6(r_s,y')=\kappa_-(r_s+V_7(y'),y')-\bar{\kappa},
\end{cases}\ee
where $R_{2}({\bf V}(r_s,y'),V_7(y'))=\sum_{i=1}^2 b_{2i}R_{0i}({\bf V}(r_s,y'),V_7(y'))$.

The following reformulation of the jump conditions \eqref{shock12} is crucial for us to solve the shock problem. We write \eqref{D23_V7_ycor} as
\be\label{shock121}\begin{cases}
F_2(y'):=\frac{1}{r_s}\p_{y_2} V_7- a_0 V_2(r_s, y')- g_2({\bf V}(r_s,y'),V_7)\equiv 0,\ \ &\forall y'\in E,\\
F_3(y'):=\p_{y_3}V_7- a_0  V_3(r_s,y')- g_3({\bf V}(r_s,y'),V_7)\equiv 0,\ \ &\forall y'\in E.
\end{cases}\ee
The following equivalent reformulation of \eqref{shock121} is crucial.
\begin{lemma}\label{equi0}
{\it Let $F_j, j=2,3$ be two $C^1$ smooth functions defined on $\overline{E}$. Then the following two statements are equivalent
\begin{enumerate}[(1)]
\item $F_2=F_3\equiv 0$ on $\overline{E}$;
\item $F_2$ and $F_3$ solve the following problem
\be\label{equi00}\begin{cases}
\frac{1}{r_s}\p_{y_2}F_3-\p_{y_3} F_2=0,\ \ &\text{in}\ E,\\
\frac{1}{r_s}\p_{y_2} F_2 + \p_{y_3} F_3=0,\ \ &\text{in}\ E,\\
F_2(\pm\theta_0, y_3)=0,\ \ &\text{on}\ \ y_3\in [-1,1],\\
F_3(y_2, \pm 1)=0,\ \ &\text{on}\ \ y_2\in [-\theta_0,\theta_0].
\end{cases}\ee
\end{enumerate}
}\end{lemma}

The first equation in \eqref{equi00} is
\be\label{shock17}
(\frac{1}{r_s}\p_{y_2} V_3-\p_{y_3}V_2)(r_s,y')=\frac{1}{a_0} (\p_{y_3}g_2-\frac{1}{r_s}\p_{y_2}g_3)({\bf V}(r_s,y'), V_7),
\ee
which yields the data on the shock for the first component of the modified vorticity.

The second equation in \eqref{equi00} gives
\begin{equation*}
\begin{split}
&(\frac{1}{r_s^2}\p_{y_2}^2 + \p_{y_3}^2) V_7(y')-a_0(\frac{1}{r_s}\p_{y_2} V_2+\p_{y_3} V_3)(r_s, y')\\
=&(\frac{1}{r_s}\p_{y_2}g_2+ \p_{y_3}g_3)({\bf V}(r_s,y'), V_7).
\end{split}
\end{equation*}
This, together with \eqref{shock400}, shows
\be\label{shock19}
\{(\frac{1}{r_s^2}\p_{y_2}^2 +\p_{y_3}^2) V_1-a_0 a_1 (\frac{1}{r_s}\p_{y_2} V_2+\p_{y_3} V_3)\}(r_s,y')=q_1({\bf V}(r_s,y'),V_7),
\ee
with
\be\no
q_1=a_1(\frac{1}{r_s}\p_{y_2}g_2+\p_{y_3}g_3)({\bf V}(r_s,y'), V_7)+(\frac{1}{r_s^2}\p_{y_2}^2 R_1+\p_{y_3}^2 R_1)({\bf V}(r_s,y'), V_7).
\ee
The condition \eqref{shock19} is used as the boundary condition on the shock front for the deformation-curl system associated with the velocity field.

The last two equations in \eqref{equi00} can be rewritten as
\be\label{shock20}\begin{cases}
(\frac{1}{r_s}\p_{y_2}V_1-a_0a_1 V_2)(r_s,\pm \theta_0,y_3)=q_2^{\pm}({\bf V}(r_s,\pm\theta_0,y_3),V_7(\pm\theta_0,y_3)),\\
(\p_{y_3} V_1-a_0 a_1 V_3)(r_s,y_2,\pm 1)= q_3^{\pm}({\bf V}(r_s,y_2,\pm 1),V_7(y_2,\pm 1)),
\end{cases}\ee
with
\be\no\begin{cases}
q_2^{\pm}({\bf V}(r_s,\cdot),V_7)(\pm\theta_0,y_3)=(\frac{1}{r_s}\p_{y_2} \{R_1({\bf V}(r_s,\cdot), V_7(\cdot))\}+ g_2({\bf V}(r_s,\cdot), V_7))(\pm \theta_0, y_3),\\
q_3^{\pm}({\bf V}(r_s,\cdot),V_7)(y_2,\pm 1)=(\p_{y_3}\{R_1({\bf V}(r_s,\cdot), V_7(\cdot))\}+g_3({\bf V}(r_s,\cdot), V_7))(y_2,\pm 1).
\end{cases}\ee
The role of \eqref{shock20} will be indicated later.

Next, we determine the modified vorticity. Rewrite \eqref{vor21} as
\be\label{vor400}
\b(D_1^{V_7}+\frac{1}{\bar{U}(D_0^{V_7})+V_1}\sum\limits_{i=2}^3 V_i D_i^{V_7}\b) \tilde{J}_1 + \mu({\bf V},V_7)\tilde{J}_1=H_0({\bf V},V_7),
\ee
where
\be\no
&&\mu({\bf V},V_7)=\displaystyle\sum_{i=2}^3 D_i^{V_7}\bigg(\frac{V_i}{\bar{U}(D_0^{V_7})+V_1}\bigg)+\frac{1}{D_0^{V_7}},\\\no
&&H_0({\bf V}, V_7)=D_3^{V_7}\b(\frac{1}{\bar{U}(D_0^{V_7})+V_1}\b)D_2^{V_7} V_5- D_2^{V_7}\b(\frac{1}{\bar{U}(D_0^{V_7})+V_1}\b)D_3^{V_7} V_5\\\no
&&\quad+D_2^{V_7}\b(\frac{\bar{B}+V_5-\frac{1}{2}(\bar{U}(D_0^{V_7})+V_1)^2-\frac{1}{2}(V_2^2+V_3^2)}{\gamma(\bar{S}+V_4)(\bar{U}(D_0^{V_7})+V_1)}\b)D_3^{V_7} V_4\\\no
&&\quad- D_3^{V_7}\b(\frac{\bar{B}+V_5-\frac{1}{2}(\bar{U}(D_0^{V_7})+V_1)^2-\frac{1}{2}(V_2^2+V_3^2)}{\gamma(\bar{S}+V_4)(\bar{U}(D_0^{V_7})+V_1)}\b)D_2^{V_7} V_4\\\no
&&\quad+D_2^{V_7}\b(\frac{(\bar{\ka}+V_6)\tilde{\rho}({\bf V})((\bar{U}(D_0^{V_7})+V_1)^2+V_2^2+V_3^2)}{(\bar{U}(D_0^{V_7})+V_1)}\b)D_3^{V_7} V_6\\\no
&&\quad- D_3^{V_7}\b(\frac{(\bar{\ka}+V_6)\tilde{\rho}({\bf V})((\bar{U}(D_0^{V_7})+V_1)^2+V_2^2+V_3^2)}{(\bar{U}(D_0^{V_7})+V_1)}\b)D_2^{V_7} V_6.
\ee
Then \eqref{shock17} gives the boundary data for $\tilde{J}_1$ at $y_1=r_s$
\be\label{vor401}
\tilde{J}_1(r_s,y')=\frac{1}{a_0}(\p_{y_3} g_2-\frac{1}{r_s}\p_{y_2}g_3)({\bf V}(r_s,y'), V_7)+g_4({\bf V}(r_s,y'), V_7(y')),
\ee
with
\be\label{def_g4}
&&g_4({\bf V}, V_7)=(D_2^{V_7}-\frac1{y_1}\p_{y_2})\{(1-\bar{\ka}^2\bar{\rho})V_3\}-(D_3^{V_7}-\p_{y_3})\{(1-\bar{\ka}^2\bar{\rho})V_2\}\\\no
&&\q-D_2^{V_7}\{((\bar{\ka}+V_6)^2\tilde{\rho}({\bf V},V_7)-\bar{\ka}^2\bar{\rho})V_3\}+D_3^{V_7}\{((\bar{\ka}+V_6)^2\tilde{\rho}({\bf V},V_7)-\bar{\ka}^2\bar{\rho})V_2\}.
\ee
On the other hand, \eqref{vor22} implies that
\be\no
&&\tilde{J}_2=D_3^{V_7} ((1-(\bar{\ka}+V_6)^2\tilde{\rho}({\bf V},V_7))V_1)\\\no
&&\q\q- D_1^{V_7}((1-(\bar{\ka}+V_6)^2\tilde{\rho}({\bf V},V_7))V_3)-\bar{U}(D_0^{V_7})D_3^{V_7}\{(\bar{\ka}+V_6)^2\tilde{\rho}({\bf V},V_7)\}\\\no
&&=\frac{V_2\tilde{J}_1+D_3^{V_7} V_5-(\bar{\ka}+V_6)\tilde{\rho}({\bf V}, V_7)((\bar{U}(D_0^{V_7})+V_1)^2+V_2^2+V_3^2)D_3^{V_7}V_6}{\bar{U}(D_0^{V_7})+V_1}\\\label{vor402}
&&\q\q-\frac{\bar{B}-\frac{1}{2}\bar{U}^2(D_0^{V_7})+V_5-\bar{U}(D_0^{V_7})V_1-\frac{1}{2}\sum_{i=1}^3 V_i^2}{\gamma (\bar{U}(D_0^{V_7})+V_1)(\bar{S}+V_4)} D_3^{V_7} V_4,\\\no
&&\tilde{J}_3=(D_1^{V_7}+\frac{1}{D_0^{V_7}})((1-(\bar{\ka}+V_6)^2\tilde{\rho}({\bf V},V_7))V_2)\\\no
&&\q\q-D_2^{V_7} ((1-(\bar{\ka}+V_6)^2\tilde{\rho}({\bf V},V_7))V_1)+\bar{U}(D_0^{V_7})D_2^{V_7}\{(\bar{\ka}+V_6)^2\tilde{\rho}({\bf V},V_7)\}\\\no
&&=\frac{V_3\tilde{J}_1-D_2^{V_7} V_5+(\bar{\ka}+V_6)\tilde{\rho}({\bf V}, V_7)((\bar{U}(D_0^{V_7})+V_1)^2+V_2^2+V_3^2)D_2^{V_7}V_6}{\bar{U}(D_0^{V_7})+V_1}\\\label{vor403}
&&\q\q+\frac{\bar{B}-\frac{1}{2}\bar{U}^2(D_0^{V_7})+V_5-\bar{U}(D_0^{V_7})V_1-\frac{1}{2}\sum_{i=1}^3 V_i^2}{\gamma (\bar{U}(D_0^{V_7})+V_1)(\bar{S}+V_4)} D_2^{V_7} V_4.
\ee

Collecting the principal terms and putting the quadratic terms on the right hand sides, one gets from direct computations and \eqref{vor402}-\eqref{vor403} that
\be\label{vor404}
&&\frac{1}{y_1}\p_{y_2} \{d_0 V_3\}-\p_{y_3} \{d_0 V_2\}=\tilde{J}_1(y)+H_1({\bf V},V_7),\\\no
&&\p_{y_3}\{d V_1+d_3 V_4\}- \p_{y_1} (d_0 V_3) = H_2({\bf V},V_7)\\\label{vor405}
&&\quad+\bar{U}(D_0^{V_7})\b\{2(\bar{\ka}+V_6)\tilde{\rho}({\bf V},V_7)D_3^{V_7}V_6+\frac{(\bar{\ka}+V_6)^2\tilde{\rho}({\bf V},V_7)}{c^2(\tilde{\rho})}D_3^{V_7}V_5\b\}\\\no
&&\quad+\frac{V_2\tilde{J}_1+D_3^{V_7} V_5-(\bar{\ka}+V_6)\tilde{\rho}({\bf V},V_7)((\bar{U}(D_0^{V_7})+V_1)^2+V_2^2+V_3^2)D_3^{V_7}V_6}{\bar{U}(D_0^{V_7})+V_1},\\\no
&&(\p_{y_1}+\frac{1}{y_1})\{d_0 V_2\}-\frac{1}{y_1} \p_{y_2} \{d V_1+d_3V_4\} =H_3({\bf V},V_7)\\\label{vor406}
&&\q-\bar{U}(D_0^{V_7})\b\{2(\bar{\ka}+V_6)\tilde{\rho}({\bf V},V_7)D_2^{V_7}V_6+\frac{(\bar{\ka}+V_6)^2\tilde{\rho}({\bf V},V_7)}{c^2(\tilde{\rho})}D_2^{V_7}V_5\b\}\\\no
&&\q+\frac{V_3\tilde{J}_1-D_2^{V_7} V_5+(\bar{\ka}+V_6)\tilde{\rho}({\bf V}, V_7)((\bar{U}(D_0^{V_7})+V_1)^2+V_2^2+V_3^2)D_2^{V_7}V_6}{\bar{U}(D_0^{V_7})+V_1},
\ee
where
\be\no
&&H_1({\bf V},V_7)=(D_3^{V_7}-\p_{y_3})\{d_0 V_2\}-(D_2^{V_7}-\p_{y_2})\{d_0 V_3\}\\\no
&&\q+ D_2^{V_7}\{((\bar{\ka}+V_6)^2\tilde{\rho}({\bf V}, V_7)-\bar{\ka}^2\bar{\rho})V_3\}-D_3^{V_7}\{((\bar{\ka}+V_6)^2\tilde{\rho}({\bf V}, V_7)-\bar{\ka}^2\bar{\rho})V_2\},\\\no
&&H_2({\bf V},V_7)=-\frac{V_5-\bar{U}(D_0^{V_7})V_1-\frac{1}{2}\sum_{j=1}^3 V_j^2}{\gamma (\bar{U}(D_0^{V_7})+V_1)(\bar{K}+V_4)} D_3^{V_7} V_4+\frac{(\bar{\ka}+V_6)^2\tilde{\rho}\bar{U}^2(D_0^{V_7})}{c^2(\tilde{\rho})}D_3^{V_7} V_1\\\no
&&\q-\bar{\ka}^2\bar{\rho}\bar{M}^2\p_{y_3}V_1-\b(\frac{\bar{B}-\frac{1}{2}\bar{U}^2(D_0^{V_7})}{\gamma (\bar{U}(D_0^{V_7})+V_1)(\bar{S}+V_4)}+\frac{(\bar{\ka}+V_6)^2\tilde{\rho}\bar{U}(D_0^{V_7})}{(\gamma-1)(\bar{S}+V_4)}\b)D_3^{V_7}V_4\\\no
&&\q+\p_{y_3} (d_3 V_4)-\frac{(\bar{\ka}+V_6)^2\tilde{\rho}\bar{U}(D_0^{V_7})}{c^2(\tilde{\rho})}\sum_{i=1}^3 V_i D_3^{V_7}V_i\\\no
&&\q+(D_1^{V_7}-\p_{y_1})\{d_0 V_3\}-(D_3^{V_7}-\p_{y_3})\{d_0 V_1\}\\\no
&&\q+ D_3^{V_7}\{((\bar{\ka}+V_6)^2\tilde{\rho}({\bf V}, V_7)-\bar{\ka}^2\bar{\rho})V_1\}-D_1^{V_7}\{((\bar{\ka}+V_6)^2\tilde{\rho}({\bf V}, V_7)-\bar{\ka}^2\bar{\rho})V_3\},\\\no
&&H_3({\bf V},V_7)=\frac{V_5-\bar{U}(D_0^{V_7})V_1-\frac{1}{2}\sum_{i=1}^3 V_i^2}{\gamma (\bar{U}(D_0^{V_7})+V_1)(\bar{S}+V_4)} D_2^{V_7} V_4\\\no
&&\quad\quad+ \b(\frac{\bar{B}-\frac{1}{2}\bar{U}^2(D_0^{V_7})}{\gamma (\bar{U}(D_0^{V_7})+V_1)(\bar{S}+V_4)}+\frac{(\bar{\ka}+V_6)^2\tilde{\rho}\bar{U}(D_0^{V_7})}{(\gamma-1)(\bar{S}+V_4)}\b)D_2^{V_7} V_4\\\no
&&\q-\frac{1}{y_1}\p_{y_2} (d_3 V_4)+\frac{(\bar{\ka}+V_6)^2\tilde{\rho}\bar{U}(D_0^{V_7})}{c^2(\tilde{\rho})}\sum_{i=1}^3 V_i D_2^{V_7}V_i+(D_2^{V_7}-\frac{1}{y_1}\p_{y_2})\{d V_1\}\\\no
&&\q-(D_1^{V_7}+\frac{1}{D_0^{V_7}}-\p_{y_1}-\frac1{y_1})\{d_0 V_2\}+(D_1^{V_7}+\frac{1}{D_0^{V_7}})\{((\bar{\ka}+V_6)^2\tilde{\rho}({\bf V}, V_7)-\bar{\ka}^2\bar{\rho})V_2\}\\\no
&&\q-D_2^{V_7}\{((\bar{\ka}+V_6)^2\tilde{\rho}({\bf V}, V_7)-\bar{\ka}^2\bar{\rho})V_1\}.
\ee
The boundary conditions on $\Sigma_2^{\pm}$ and $\Sigma_3^{\pm}$, \eqref{slip1}, become
\be\label{slip2}\begin{cases}
V_2(y_1,\pm\theta_0, y_3)=0,\ \ &\text{on }\Sigma_2^{\pm},\\
V_3(y_1,y_2,\pm 1)=0,\ \ &\text{on }\Sigma_3^{\pm}.
\end{cases}\ee

Furthermore, the equation \eqref{den15} can be rewritten as
\be\label{den20}
d_1\p_{y_1} V_1+ \frac{1}{y_1} \p_{y_2} V_2 + \p_{y_3} V_3 +\frac{V_1}{y_1}+ d_2 V_1=d_5 V_5+ G_0({\bf V},V_7),
\ee
with
\be\no
&&G_0({\bf V},V_7)=\mathbb{F}({\bf V},V_7)-\b(d_1(D_0^{V_7})D_1^{V_7} V_1-d_1(y_1)\p_{y_1} V_1\b)-(D_2^{V_7} V_2-\frac{1}{y_1}\p_{y_2}V_2)\\\no
&&\quad-(D_3^{V_7}V_3-\p_{y_3} V_3)- \b((\frac{1}{D_0^{V_7}} +d_2(D_0^{V_7}))V_1-(\frac{1}{y_1}+d_2(y_1)) V_1\b),
\ee
and
\be\no
&&\bar{c}^2(D_0^{V_7})\mathbb{F}({\bf V},V_7)=-(\gamma-1)(D_1^{V_7} V_1+\frac{V_1}{D_0^{V_7}}) V_5\\\no
&&\quad\quad+ (\bar{U}'(D_0^{V_7})+D_1^{V_7} V_1)\b(\frac{\gamma+1}{2} V_1^2+\frac{\gamma-1}{2}(V_2^2+V_3^2)\b)\\\no
&&\quad\quad+\frac{(\gamma-1)(\bar{U}(D_0^{V_7})+V_1)}{2 D_0^{V_7}}\sum_{i=1}^3 V_i^2+\bar{U}(D_0^{V_7})\{(\gamma+1)V_1 D_1^{V_7} V_1+\frac{\gamma-1}{D_0^{V_7}}V_1^2\}\\\no
&&\quad\quad- (\gamma-1)\b(V_5-\frac{1}{2}\sum_{i=1}^3 V_i^2-\bar{U}(D_0^{V_7})V_1\b)(D_2^{V_7} V_2+D_3^{V_7} V_3)\\\no
&&\quad\quad
+(V_2^2 D_2^{V_7} V_2 +V_3^2 D_3^{V_7}V_3)+ (\bar{U}(D_0^{V_7})+V_1)(V_2 D_1^{V_7} V_2+V_3 D_1^{V_7} V_3)\\\no
&&\quad\quad+V_2((\bar{U}(D_0^{V_7})+V_1)D_2^{V_7} V_1 +V_3 D_2^{V_7} V_3)\\\no
&&\quad\quad+ V_3((\bar{U}(D_0^{V_7})+V_1)D_3^{V_7} V_1 +V_2 D_3^{V_7} V_2).
\ee

Finally, the boundary condition \eqref{pres2} at the exit becomes
\be\no
&&(d(r_2)V_1+d_3(r_2)V_4)(r_2,y')=-\frac{\epsilon T_e(y')}{(\bar{\rho}\bar{U})(r_2)}+\frac{(1+\bar{\kappa}^2 \bar{\rho}\bar{M}^2)W_{5}(r_2,y')}{\bar{U}(r_2)}\\\label{pres3}
&&\quad\quad-\bar{\kappa}(\bar{\rho}\bar{U})(r_2) V_6(r_2,\cdot)-\frac{E({\bf V}(r_2,y'))}{(\bar{\rho}\bar{U})(r_2)},
\ee
where
\be\no
&&E({\bf V})(r_2,y')=\tilde{P}({\bf V},V_7)+ \frac12 (\bar{\kappa}+V_{6})^2(\tilde{\rho}({\bf V},V_7))^2((\bar{U}+V_{1})^2 + V_{2}^2 +V_{3}^2)(r_2,y')\\\no
&&\quad-(\bar{P}+\frac12 \bar{\kappa}^2 \bar{\rho}^2 \bar{U}^2)(r_2)+ (\bar{\rho}\bar{U})(r_2)(d(r_2) V_1 + d_{3}(r_2) V_{4})(r_2,y')\\\label{err3}
&&\quad -\bar{\rho}(r_2)(1+\bar{\kappa}^2 \bar{\rho}\bar{M}^2(r_2))V_5(r_2,y')-\bar{\kappa}(\bar{\rho}\bar{U})(r_2) V_6(r_2,y').
\ee

Therefore, after the coordinates transformation \eqref{coor}, the transonic shock problem \eqref{mhd-cyl} with \eqref{super1}-\eqref{pressure} and \eqref{rh} is equivalent to solving the following problem.

{\bf Problem S.} Find a function $V_7$ defined on $E$ and vector functions $(V_1,\cdots, V_6)$ defined on the $\mathbb{D}$, which solve the transport equations \eqref{ber31}, \eqref{vor400},\eqref{vor404}-\eqref{vor406} and \eqref{den20} with boundary conditions \eqref{shock400}, \eqref{shock19}-\eqref{shock20}, \eqref{vor401},\eqref{slip2} and \eqref{pres3}.

Theorem \ref{existence} then follows directly from the following result.
\begin{theorem}\label{main}
{\it Assume that the compatibility conditions \eqref{pressure-cp} and \eqref{super3} hold. There exists a small constant $\epsilon_0>0$ depending only on the background solution and the boundary data such that if $0\leq \epsilon<\epsilon_0$, the problem \eqref{ber31},\eqref{vor400},\eqref{vor404}-\eqref{vor406},\eqref{den20} with boundary conditions \eqref{shock400}, \eqref{shock19}-\eqref{shock20},\eqref{vor401},\eqref{slip2} and \eqref{pres3} has a unique solution $(V_1,V_2,V_3,V_4,V_5,V_6)(y)$ with the shock front $\mlS: y_1=V_7(y')$ satisfying the following properties.
\begin{enumerate}[(1)]
\item The function $V_7(y')\in C^{3,\alpha}(\overline{E})$ satisfies
\be\no
\|V_7(y')\|_{C^{3,\alpha}(\overline{E})}\leq C_*\epsilon,
\ee
and
\be\label{shock91}\begin{cases}
\p_{y_2}V_7(\pm\theta_0, y_3)=\p_{y_2}^3V_7(\pm\theta_0, y_3)=0,\ \ &\forall y_3\in [-1,1],\\
\p_{y_3}V_7(y_2,\pm 1)=\p_{y_3}^3V_7(y_2, \pm 1)=0,\ \ &\forall y_2\in [-\theta_0,\theta_0],
\end{cases}\ee
where $C_*$ depends only on the background solution, the supersonic incoming flow, and the exit pressure.
\item The solution $(V_1,V_2,V_3,V_4,V_5,V_6)(y)\in C^{2,\alpha}(\overline{\mathbb{D}})$ satisfies
\be\no
\sum_{j=1}^6\|V_j\|_{C^{2,\alpha}(\overline{\mathbb{D}})}\leq C_*\epsilon,
\ee
and the compatibility conditions
\be\label{sub52}
\begin{cases}
(V_2,\p_{y_2}^2V_2)(y_1,\pm\theta_0,y_3)=\p_{y_2}(V_1,V_3,V_4,V_5,V_6))(y_1,\pm\theta_0,y_3)=0, &\text{on }\Sigma_2^{\pm},\\
(V_3,\p_{y_3}^2V_3)(y_1,y_2,\pm 1)=\p_{y_3}(V_1,V_2,V_4,V_5,V_6))(y_1,y_2,\pm 1)=0,\ &\text{on }\Sigma_3^{\pm}.
\end{cases}\ee
\end{enumerate}
}\end{theorem}

\section{Proof of Theorem \ref{main}}\label{proof}

We proceed to prove Theorem \ref{main}. The solution class $\Xi$ consists of the vector functions $(V_1,\cdots, V_6,V_7)\in (C^{2,\alpha}(\overline{\mathbb{D}}))^6\times C^{3,\alpha}(\overline{E})$ satisfying the estimate
\be\no
\|({\bf V},V_7)\|_{\Xi}:=\sum_{i=1}^6 \|V_i\|_{C^{2,\alpha}(\overline{\mathbb{D}})}+\|V_7\|_{C^{3,\alpha}(\overline{E})}\leq \delta_0,
\ee
and the following compatibility conditions (which is precisely \eqref{shock91} and \eqref{sub52})
\be\label{class2}\begin{cases}
(V_2,\p_{y_2}^2V_2,\p_{y_2}(V_1,V_3,V_4,V_5,V_6))(y_1,\pm\theta_0,y_3)=0,\ \ &\text{on }\Sigma_2^{\pm},\\
(V_3,\p_{y_3}^2V_3,\p_{y_3}(V_1,V_2,V_4,V_5,V_6))(y_1,y_2,\pm 1)=0,\ \ &\text{on }\Sigma_3^{\pm},\\
(\p_{y_2}V_7,\p_{y_2}^3V_7)(\pm\theta_0,y_3)=0,\ \ &\text{on }y_3\in [-1,1],\\
(\p_{y_3}V_7,\p_{y_3}^3V_7)(y_2,\pm 1)=0,\ \  &\text{on } y_2\in [-\theta_0,\theta_0],
\end{cases}\ee
with $\delta_0$ being a suitably small positive constant to be determined later.

For any given $(\hat{{\bf V}},\hat{V}_7)\in \Xi$, we will define an operator $\mathcal{T}$ mapping $\Xi$ to itself, and the unique fixed point of $\mathcal{T}$ will solve the {\bf Problem S}.

{\bf Step 1.} The shock front is uniquely determined by the following algebraic equation:
\be\label{shock41}
V_7(y')=a_1^{-1} V_1(r_s,y')-a_1^{-1} R_1(\hat{{\bf V}}(r_s,y'),\hat{V}_7),
\ee
provided that $V_1(r_s,y')$ is obtained.

{\bf Step 2.} We solve the transport equations for the Bernoulli quantity and the entropy, respectively. The Bernoulli's quantity $V_5$ and the function $V_6$ will be determined by (See \eqref{ber31})
\be\lab{ber41}\begin{cases}
\bigg(D_1^{\hat{V}_7}+\frac{\hat{V}_2}{\bar{U}(D_0^{\hat{V}_7})+\hat{V}_1}D_2^{\hat{V}_7}+\frac{\hat{V}_3}{\bar{U}(D_0^{\hat{V}_7})+\hat{V}_1}D_3^{\hat{V}_7}\bigg)(V_5,V_6)=0,\\
V_5(r_s,y_2,y_3)=B_-(r_s+\hat{V}_7(y'),y')-\bar{B}_-,\\
V_6(r_s,y_2,y_3)=\kappa_-(r_s+\hat{V}_7(y'),y')-\bar{\kappa}_-.
\end{cases}\ee

Set
\be\label{char}\begin{cases}
K_2(y):=\frac{r_2-r_s-\hat{V}_7}{r_2-r_s}\frac{\hat{V}_2}{D_0^{\hat{V}_7}(\bar{U}(D_0^{\hat{V}_7})+\hat{V}_1)
+\frac{y_1-r_2}{r_2-r_s}(\hat{V}_2\p_{y_2}\hat{V}_7+D_0^{\hat{V}_7}\hat{V}_3\p_{y_3}\hat{V}_7)},\\
K_3(y):=\frac{r_2-r_s-\hat{V}_7}{r_2-r_s}\frac{D_0^{\hat{V}_7}\hat{V}_3}{\bar{U}(D_0^{\hat{V}_7})+\hat{V}_1
+\frac{y_1-r_2}{r_2-r_s}(\hat{V}_2\p_{y_2}\hat{V}_7+D_0^{\hat{V}_7}\hat{V}_3\p_{y_3}\hat{V}_7)}.
\end{cases}\ee
Then $K_2,K_3\in C^{2,\alpha}(\overline{\mathbb{D}})$ for any $(\hat{{\bf V}},\hat{V}_7)\in \Xi$. Define the trajectory by solving the ODE system
\begin{eqnarray}\label{char1} \left\{\begin{array}{ll}
\frac{d \bar{y}_2(\tau; y)}{d\tau}=K_2(\tau,\bar{y}_2(\tau;y),\bar{y}_3(\tau;y)),\ \ \forall \tau\in [r_s,r_2],\\
\frac{d \bar{y}_3(\tau; y)}{d\tau}=K_3(\tau,\bar{y}_2(\tau;y),\bar{y}_3(\tau;y)),\ \ \forall \tau\in [r_s,r_2],\\
\bar{y}_2(y_1; y)=y_2,\ \bar{y}_3(y_1;y)=y_3.
\end{array}\right. \end{eqnarray}
Denote $\vec{\beta}(y)=(\beta_2(y),\beta_3(y))=(\bar{y}_2(r_s;y),\bar{y}_3(r_s;y))$. Since $(\hat{{\bf V}},\hat{V}_7)\in \Xi$ satisfies the compatibility conditions \eqref{class2}, then
\be\label{char201}\begin{cases}
K_2(y_1,\pm\theta_0,y_3)=\p_{y_2} K_3(y_1,\pm\theta_0,y_3)=0,\ \ &\text{on }\Sigma_2^{\pm},\\
K_3(y_1,y_2,\pm 1)=\p_{y_3} K_2(y_1,y_2, \pm 1)=0,\ \ &\text{on }\Sigma_3^{\pm}.
\end{cases}\ee
According to the uniqueness of the solution to \eqref{char1} and \eqref{char201}, there hold
\be\label{char20}\begin{cases}
\bar{y}_2(\tau;y_1,\pm\theta_0,y_3)=\pm\theta_0,\ \ &\forall \tau\in [r_s,r_2], (y_1,y_3)\in \Sigma_2^{\pm},\\
\bar{y}_3(\tau;y_1,y_2,\pm 1)=\pm 1,\ \ &\forall \tau\in [r_s,r_2], (y_1,y_2)\in \Sigma_3^{\pm},
\end{cases}\ee
and
\be\label{char2}\begin{cases}
\beta_2(y_1,\pm\theta_0,y_3)=\pm\theta_0,\ &\forall  (y_1,y_3)\in \Sigma_2^{\pm},\\
\beta_3(y_1,y_2,\pm 1)=\pm 1,\ &\forall  (y_1,y_2)\in \Sigma_3^{\pm}.
\end{cases}\ee
The existence and uniqueness of $(\bar{y}_2(\tau;y),\bar{y}_3(\tau;y))$ on the whole interval $[r_s,r_2]$ follow from the standard theory of systems of ordinary differential equations and \eqref{char20}. Furthermore, it holds
\be\no
\sum_{j=2}^3\|\beta_j(y)-y_j\|_{C^{2,\alpha}(\overline{\mathbb{D}})}\leq C_*\|(\hat{{\bf V}},\hat{V}_7)\|_{\Xi}.
\ee 

Furthermore, the functions $\beta_2,\beta_3$ satisfy the conditions:
\be\label{char3}\begin{cases}
\p_{y_2}\beta_3(y_1,\pm\theta_0, y_3)=0,\ \ &\text{on }\Sigma_2^{\pm},\\
\p_{y_3}\beta_2(y_1,y_2,\pm 1)=0,\ \ &\text{on }\Sigma_3^{\pm}.
\end{cases}\ee



Since $V_5$ and $V_6$ are conserved along the trajectory, one has
\be\no
&&V_5(y)=V_5(r_s,\vec{\beta}(y))=B_-(r_s+\hat{V}_7(\vec{\beta}(y)),\vec{\beta}(y))-\bar{B}_-,\\\no
&&V_6(y)=V_6(r_s,\vec{\beta}(y))=\kappa_-(r_s+\hat{V}_7(\vec{\beta}(y)),\vec{\beta}(y))-\bar{\kappa}_-.
\ee
Thus $V_5$ and $V_6$ can be regarded as high order terms with the following estimate
\be\label{ber43}
&&\sum_{i=5}^6\|V_i\|_{C^{2,\alpha}(\overline{\mathbb{D}})}\leq C_*\epsilon (\|\hat{V}_7\|_{C^{2,\alpha}(\overline{E})}+\sum_{j=2}^3\|\beta_j\|_{C^{2,\alpha}(\overline{\mathbb{D}})})\\\no
&&\leq C_*(\epsilon+\epsilon\|(\hat{{\bf V}}, \hat{V}_7)\|_{\Xi})\leq C_*(\epsilon+\epsilon \delta_0).
\ee
It follows from \eqref{super5},\eqref{class2} and \eqref{char3} that the following compatibility conditions hold
\be\label{ber44}\begin{cases}
\p_{y_2}(V_5,V_6)(y_1,\pm \theta_0, y_3)=0,\ \ \text{on }\Sigma_2^{\pm},\\
\p_{y_3}(V_5,V_6)(y_1,y_2,\pm 1)=0,\ \ \text{on }\Sigma_3^{\pm}.
\end{cases}\ee

The function $V_4$ satisfies
\be\no\begin{cases}
\b(D_1^{\hat{V}_7}+\frac{\hat{V}_2}{\bar{U}(D_0^{\hat{V}_7})+\hat{V}_1}D_2^{\hat{V}_7}+\frac{\hat{V}_3}{\bar{U}(D_0^{\hat{V}_7})+\hat{V}_1}D_3^{\hat{V}_7}\b)V_4=0,\\
V_4(r_s,y')= a_2 V_7(y')+R_2(\hat{{\bf V}}(r_s,y'),\hat{V}_7(y')).
\end{cases}\ee

By the characteristic method and the equation \eqref{shock41}, one has
\be\label{ent42}
&&V_4(y)=V_4(r_s,\vec{\beta}(y))\\\no
&&=a_2 V_7(\vec{\beta}(y))+ R_2(\hat{{\bf V}}(r_s,\vec{\beta}(y)),\hat{V}_7(\vec{\beta}(y)))\\\no
&&=a_2 V_7(y')+ a_2(V_7(\vec{\beta}(y))-V_7(y'))+ R_2(\hat{{\bf V}}(r_s,\vec{\beta}(y)),\hat{V}_7(\vec{\beta}(y)))\\\no
&&=\frac{a_2}{a_1} V_1(r_s,y')+a_2(V_7(\vec{\beta}(y))-V_7(y'))+R_3(\hat{{\bf V}}(r_s,\vec{\beta}(y)),\hat{V}_7(\vec{\beta}(y))),
\ee
and
\be\no
R_3(\hat{{\bf V}}(r_s,y'),\hat{V}_7)=R_2(\hat{{\bf V}}(r_s,y'),\hat{V}_7)-\frac{a_2 }{a_1}R_1(\hat{{\bf V}}(r_s,y'),\hat{V}_7).
\ee

Since $V_7(y')$ is still unknown, one may rewrite \eqref{ent42} as
\be\label{ent43}
V_4(y_1,y')=\frac{a_2}{a_1} V_1(r_s,y')+R_4(\hat{{\bf V}}(r_s,\vec{\beta}(y)),\hat{V}_7(\vec{\beta}(y))),
\ee
with
\be\no
R_4=a_2(\hat{V}_7(\vec{\beta}(y))-\hat{V}_7(y'))+R_3(\hat{{\bf V}}(r_s,\vec{\beta}(y)),\hat{V}_7(\vec{\beta}(y))).
\ee
Therefore $V_4$ is decomposed as a scalar multiple of $V_1(r_s,y')$ with high order terms satisfying
\be\no
&&\|V_4\|_{C^{2,\alpha}(\overline{\mathbb{D}})}\leq C_*\|V_1(r_s,\cdot)\|_{C^{2,\alpha}(\overline{E})}+\|R_4\|_{C^{2,\alpha}(\overline{\mathbb{D}})}\\\no
&&\leq C_*(\|V_1(r_s,\cdot)\|_{C^{2,\alpha}(\overline{E})}+\|\hat{V}_7|_{C^{3,\alpha}(\overline{E})}\sum_{j=2}^3\|\beta_j(y)-y_j\|_{C^{2,\alpha}(\overline{\mathbb{D}})})\\\no
&&\q+ C_*(\epsilon \|(\hat{{\bf V}}, \hat{V}_7)\|_{\Xi}+\|(\hat{{\bf V}}, \hat{V}_7)\|_{\Xi}^2)\leq C_*\|V_1(r_s,\cdot)\|_{C^{2,\alpha}(\overline{E})}+ C_*(\epsilon \delta_0+ \delta_0^2).
\ee
Furthermore, since $(\hat{{\bf V}}, \hat{V}_7)\in\Xi$ satisfies the compatibility conditions \eqref{class2} and the upcoming supersonic flow satisfies \eqref{super5}, using the expression of $f_2, f_3, f, R_{0i}, i=1,2$, one could verify by direct but tedious computations that
\be\label{j211}\begin{cases}
f_2(\hat{{\bf V}}(r_s,y'),\hat{V}_7)\}|_{y_2=\pm\theta_0}=\p_{y_2}^2\{f_2(\hat{{\bf V}}(r_s,y'),\hat{V}_7)\}|_{y_2=\pm\theta_0}=0,\ \ \\
\p_{y_2}\{f_3(\hat{{\bf V}}(r_s,y'),\hat{V}_7)\}|_{y_2=\pm\theta_0}=  \p_{y_2}\{f(\hat{{\bf V}}(r_s,y'),\hat{V}_7)\}|_{y_2=\pm\theta_0}=0,\ \\
f_3(\hat{{\bf V}}(r_s,y'),\hat{V}_7)\}|_{y_1=\pm 1}=\p_{y_3}^2\{f_3(\hat{{\bf V}}(r_s,y'),\hat{V}_7)\}|_{y_3=\pm 1}=0,\ \ \\
\p_{y_3}\{f_2(\hat{{\bf V}}(r_s,y'),\hat{V}_7)\}|_{y_3=\pm 1}=  \p_{y_3}\{f(\hat{{\bf V}}(r_s,y'),\hat{V}_7)\}|_{y_3=\pm 1}=0,\
\end{cases}\ee
and for all $j=1,2$
\be\label{rv11}\begin{cases}
\p_{y_2}\{ R_{0j}(\hat{{\bf V}}(r_s,y'),\hat{V}_7)\}|_{y_2=\pm\theta_0}=0,\ \ &\forall y_3\in [-1,1],\\
\p_{y_3}\{ R_{0j}(\hat{{\bf V}}(r_s,y'),\hat{V}_7)\}|_{y_3=\pm 1}=0,\ \ &\forall y_2\in [-\theta_0,\theta_0].
\end{cases}\ee
Thus for $k=1,2$
\be\label{rv41}\begin{cases}
\p_{y_2}\{ R_{k}(\hat{{\bf V}}(r_s,y'),\hat{V}_7)\}|_{y_2=\pm\theta_0}=0,\ &\forall y_3\in [-1,1],\\
\p_{y_3}\{ R_{k}(\hat{{\bf V}}(r_s,y'),\hat{V}_7)\}|_{y_3=\pm 1}=0,\ &\forall y_2\in [-\theta_0,\theta_0].
\end{cases}\ee
These, together with \eqref{char2} and \eqref{char3}, imply that
\be\label{rv42}
\begin{cases}
\p_{y_2}\{R_4(\hat{{\bf V}}(r_s,\vec{\beta}(y)),\hat{V}_7(\vec{\beta}(y)))\}(y_1,\pm \theta_0, y_3)=0,\ &\text{on }\Sigma_2^{\pm},\\
\p_{y_3}\{R_4(\hat{{\bf V}}(r_s,\vec{\beta}(y)),\hat{V}_7(\vec{\beta}(y)))\}(y_1,y_2,\pm 1)=0,\ &\text{on }\Sigma_3^{\pm},
\end{cases}\ee
and
\be\no\begin{cases}
\p_{y_2} V_4(y_1,\pm \theta_0, y_3)=\frac{a_2}{a_1}\p_{y_2} V_1(r_s,\pm\theta_0,y_3),\ &\text{on }\Sigma_2^{\pm},\\
\p_{y_3} V_4(y_1,y_2,\pm 1)=\frac{a_2}{a_1}\p_{y_3} V_1(r_s,y_2,\pm 1),\ &\text{on }\Sigma_3^{\pm}.
\end{cases}\ee


{\bf Step 3.} We solve the transport equation for the first component of the vorticity. Due to \eqref{vor400} and \eqref{vor401}, it suffices to consider the following problem:
\be\label{vor501}\begin{cases}
\bigg(D_1^{\hat{V}_7}+\displaystyle\sum_{i=2}^3\frac{\hat{V}_i D_i^{\hat{V}_7}}{\bar{U}(D_0^{\hat{V}_7})+\hat{V}_1}\bigg) \tilde{J}_1 + \mu(\hat{{\bf V}},\hat{V}_7)\tilde{J}_1=H_0(\hat{{\bf V}},\hat{V}_7),\\
\tilde{J}_1(r_s,y')=R_6(\hat{{\bf V}}(r_s,y'), \hat{V}_7(y')),
\end{cases}\ee
where
\be\no
&&R_6(\hat{{\bf V}}(r_s,y'), \hat{V}_7(y'))=\frac{1}{a_0}(\p_{y_3} \{g_2(\hat{{\bf V}}(r_s,y'), \hat{V}_7(y'))\}-\frac{1}{r_s}\p_{y_2}\{g_3(\hat{{\bf V}}(r_s,y'), \hat{V}_7(y'))\})\\\no
&&\quad\quad\quad+g_4(\hat{{\bf V}}(r_s,y'), \hat{V}_7(y')).
\ee
Since $(\hat{{\bf V}},\hat{V}_7)\in \Xi$ satisfies the compatibility conditions \eqref{class2}, using the first formula in \eqref{g21}, \eqref{g31} and \eqref{def_g4}, one can verify that
\be\label{vor5011}\begin{cases}
\tilde{J}_1(r_s,\pm \theta_0, y_3)=0, & \forall y_3\in [-1,1],\\
\tilde{J}_1(r_s,y_2,\pm 1)=0,& \forall y_2\in [-\theta_0,\theta_0],\\
H_0(\hat{{\bf V}},\hat{V}_7)(y_1,\pm\theta_0, y_3)=0,\ \ &\text{on }\Sigma_2^{\pm},\\
H_0(\hat{{\bf V}},\hat{V}_7)(y_1,y_2,\pm 1)=0,\ \ &\text{on }\Sigma_3^{\pm}.
\end{cases}\ee

Integrating the equation in \eqref{vor501} along the trajectory $(\tau,\bar{y}_2(\tau;y),\bar{y}_3(\tau;y))$ yields
\begin{eqnarray}\label{vor502}
&&\tilde{J}_1(y)= R_6(\vec{\beta}(y)) e^{-\int_{r_s}^{y_1} \mu(\hat{{\bf V}},\hat{V}_7)(t;\bar{y}_2(t;y),\bar{y}_3(t;y))dt}\\\nonumber
&&\quad\quad\quad+ \int_{r_s}^{y_1} H_0(\hat{{\bf V}},\hat{V}_7)(\tau,\bar{y}_2(\tau;y),\bar{y}_3(\tau;y)) e^{-\int_{\tau}^{y_1} \mu(\hat{{\bf V}},\hat{V}_7)(t;\bar{y}_2(t;y),\bar{y}_3(t;y))dt} d\tau.
\end{eqnarray}
Thus the following estimate holds
\be\no
&&\|\tilde{J}_1\|_{C^{1,\alpha}(\overline{\mathbb{D}})}\leq C_*(\|\tilde{\omega}_1(r_s,\cdot)\|_{C^{1,\alpha}(\overline{E})}+\|H_0(\hat{{\bf V}},\hat{V}_7)\|_{C^{1,\alpha}(\overline{\mathbb{D}})})\\\no
&&\leq C_*(\epsilon\|(\hat{{\bf V}}, \hat{V}_7)\|_{\Xi}+\|(\hat{{\bf V}}, \hat{V}_7)\|_{\Xi}^2)\leq C_*(\epsilon\delta_0+\delta_0^2).
\ee
Also \eqref{char2},\eqref{char3}, \eqref{vor5011} and \eqref{vor502} imply the following compatibility conditions
\be\label{vor5041}\begin{cases}
\tilde{J}_1(y_1,\pm\theta_0, y_3)=0,\ \ &\text{on }\Sigma_2^{\pm},\\
\tilde{J}_1(y_1,y_2, \pm 1)=0,\ \ &\text{on }\Sigma_3^{\pm}.
\end{cases}\ee

Substituting \eqref{vor502} and \eqref{ent43} into \eqref{vor404}-\eqref{vor406} yields
\be\label{vor504}
&&\frac{1}{y_1}\p_{y_2} \{d_0 V_3\}-\p_{y_3}\{d_0 V_2\}=G_1(\hat{{\bf V}},\hat{V}_7),\\\label{vor505}
&&\p_{y_3} \{d V_1+\frac{a_2}{a_1}d_3 V_1(r_s,y')\}- \p_{y_1} \{d_0 V_3\}=G_2(V_5,V_6;\hat{{\bf V}},\hat{V}_7),\\\label{vor506}
&&(\p_{y_1}+\frac{1}{y_1})\{d_0 V_2\}-\frac{1}{y_1} \p_{y_2} \{d V_1+\frac{a_2}{a_1}d_3V_1(r_s,y')\}\no\\
&&\q\q\q=G_3(V_5,V_6;\hat{{\bf V}},\hat{V}_7),
\ee
where
\be\no
&&G_1(\hat{{\bf V}},\hat{V}_7)=\tilde{J}_1(y)+H_1(\hat{{\bf V}},\hat{V}_7),\\\no
&&G_2(V_5,V_6;\hat{{\bf V}},\hat{V}_7)=H_2(\hat{{\bf V}},\hat{V}_7)+ d_3(y_1) \p_{y_3} \{R_4(\hat{{\bf V}},\hat{V}_7)\},\\\no
&&\quad+\bar{U}(D_0^{\hV_7})\b\{2(\bar{\ka}+V_6)\tilde{\rho}({\bf \hV},\hV_7)D_3^{\hV_7}V_6+\frac{(\bar{\ka}+V_6)^2\tilde{\rho}({\bf \hV},\hV_7)}{c^2(\tilde{\rho})}D_3^{\hV_7}V_5\b\}\\\no
&&\quad+\frac{\hV_2\tilde{J}_1+D_3^{\hV_7} V_5-(\bar{\ka}+V_6)\tilde{\rho}({\bf \hV},\hV_7)((\bar{U}(D_0^{\hV_7})+\hV_1)^2+\hV_2^2+\hV_3^2)D_3^{\hV_7}V_6}{\bar{U}(D_0^{\hV_7})+\hV_1},\\\no
&&G_3(V_5,V_6;\hat{{\bf V}},\hat{V}_7)=H_3(\hat{{\bf V}},\hat{V}_7)-d_3(y_1)  \{R_4(\hat{{\bf V}},\hat{V}_7)\}\\
&&\q-\bar{U}(D_0^{\hV_7})\b\{2(\bar{\ka}+V_6)\tilde{\rho}({\bf \hV},\hV_7)D_2^{\hV_7}V_6+\frac{(\bar{\ka}+V_6)^2\tilde{\rho}({\bf \hV},\hV_7)}{c^2(\tilde{\rho})}D_2^{\hV_7}V_5\b\}\\\no
&&\q+\frac{\hV_3\tilde{J}_1-D_2^{\hV_7} V_5+(\bar{\ka}+V_6)\tilde{\rho}({\bf \hV}, \hV_7)((\bar{U}(D_0^{\hV_7})+\hV_1)^2+\hV_2^2+\hV_3^2)D_2^{\hV_7}V_6}{\bar{U}(D_0^{\hV_7})+\hV_1}.
\ee

Using \eqref{class2}, \eqref{ber44}, \eqref{vor5041} and \eqref{rv42}, one can further verify the following compatibility conditions:
\be\label{cp20}\begin{cases}
(G_1(\hat{{\bf V}},\hat{V}_7), G_3(V_5,V_6;\hat{{\bf V}},\hat{V}_7),\p_{y_2} G_2(V_5,V_6;\hat{{\bf V}},\hat{V}_7))|_{y_2=\pm\theta_0}=0,\\
(G_1(\hat{{\bf V}},\hat{V}_7), G_2(V_5,V_6;\hat{{\bf V}},\hat{V}_7),\p_{y_3}G_3(V_5,V_6;\hat{{\bf V}},\hat{V}_7))|_{y_3=\pm 1}=0,
\end{cases}\ee

Furthermore, \eqref{den20} implies that
\be\label{den251}
\q\q \q d_1\p_{y_1} V_1+\frac{1}{y_1}\p_{y_2} V_2 +\p_{y_3} V_3+ \frac{V_1}{y_1}+ d_2 V_1=d_5 V_5+G_0(\hat{{\bf V}},\hat{V}_7).
\ee

It follows from \eqref{ent43} and \eqref{pres3} that
\be\label{pres4}
d(r_2)V_1(r_2,y')+ \frac{a_2}{a_1}d_{3}(r_2) V_1(r_2,y')=q_4(y'),
\ee
where
\be\no
&&q_4(y')=-\frac{a_2}{a_1}d_{3}(r_2) R_4(\hat{{\bf V}}(r_s,\vec{\beta}(r_2,y')),\hat{V}_7(\vec{\beta}(r_2,y')))-\frac{\epsilon T_e(y')}{(\bar{\rho}\bar{U})(r_2)}\\\no
&&\quad+\frac{(1+\bar{\kappa}^2 \bar{\rho}\bar{M}^2)\hat{V}_{5}(r_2,y')}{\bar{U}(r_2)}-\bar{\kappa}(\bar{\rho}\bar{U})(r_2) \hat{V}_6(r_2,\cdot)-\frac{E(\hat{{\bf V}}(r_2,y'))}{(\bar{\rho}\bar{U})(r_2)}.
\ee
And using \eqref{pressure-cp} and the explicit expression of $E(\hat{{\bf V}}(r_2,y'))$ in \eqref{err3}, one has
\be\no\begin{cases}
\p_{y_2} q_4(\pm\theta_0, y_3)=0,\ \ &\forall y_3\in [-1,1],\\
\p_{y_3} q_4(y_2, \pm 1)=0,\ \ &\forall y_2\in [-\theta_0,\theta_0].
\end{cases}\ee

{\bf Step 4.} We have derived a deformation-curl system for the velocity field which consists of the equations \eqref{den251}, \eqref{vor504}-\eqref{vor506} supplemented with the boundary conditions \eqref{shock19}, \eqref{slip2}, \eqref{pres4} and \eqref{shock20}, where $q_1$ and $q_i^{\pm} (i=2,3)$ are evaluated at $(\hat{{\bf V}},\hat{V}_7)$. However, due to the linearization, the vector field $(G_1,G_2,G_3)(\hat{{\bf V}},\hat{V}_7)$ may not be divergence free and thus the solvability condition of the curl system \eqref{vor504}-\eqref{vor506} does not hold in general. To overcome this obstacle, we first consider the following enlarged deformation-curl system, which includes an additional new unknown function $\Pi$ with homogeneous Dirichlet boundary conditions for $\Pi$:
\be\label{den32}\q\begin{cases}
d_1\p_{y_1} V_1+\frac{1}{y_1}\p_{y_2} V_2 +\p_{y_3} V_3+ \frac{V_1}{y_1}+ d_2 V_1=d_5 V_5+G_0(\hat{{\bf V}},\hat{V}_7),\\
\frac{1}{y_1}\p_{y_2}\{d_0 V_3\}- \p_{y_3} \{d_0V_2\}+\p_{y_1} \Pi=G_1(\hat{{\bf V}},\hat{V}_7),\\
\p_{y_3} \{d V_1+\frac{a_2}{a_1} d_3(y_1) V_1(r_s,y')\}-\p_{y_1} \{d_0V_3\}+ \frac1{y_1}\p_{y_2}\Pi= G_2(V_5,V_6;\hat{{\bf V}},\hat{V}_7),\\
(\p_{y_1}+\frac{1}{y_1})\{d_0 V_2\}-\frac{1}{y_1}\p_{y_2} \{d V_1+ \frac{a_2}{a_1}d_3(y_1) V_1(r_s,y')\}+\p_{y_3}\Pi = G_3(V_5,V_6;\hat{{\bf V}},\hat{V}_7),\
\end{cases}\ee
and
\be\label{den321}\begin{cases}
V_2(y_1,\pm\theta_0,y_3)=\Pi(y_1,\pm \theta_0, y_3)=0,\ &\text{on }\Sigma_2^{\pm},\\
V_3(y_1,y_2,\pm 1)=\Pi(y_1,y_2,\pm 1)= 0,\ &\text{on }\Sigma_3^{\pm}, \\
\Pi(r_s,y')=\Pi(r_2,y')=0,\ &\forall y'\in E,\\
d(r_2)V_1(r_2,y')+\frac{a_2}{a_1} d_{3}(r_2) V_1(r_s,y')=q_4(y'),\ &\forall y'\in E,\\
\{(\frac{1}{r_s^2}\p_{y_2}^2+\p_{y_3}^2) V_1-a_0 a_1 (\frac{1}{r_s}\p_{y_2} V_2+\p_{y_3}V_3)\}(r_s,y')= q_1(\hat{{\bf V}},\hat{V}_7),&\forall y'\in E\\
\b(\frac{1}{r_s}\p_{y_2}V_1-a_0a_1 V_2\b)(r_s,\pm \theta_0,y_3)=0,\ &\forall y_3\in [-1,1],\\
(\p_{y_3} V_1-a_0 a_1 V_3)(r_s,y_2,\pm 1)= 0,\  &\forall y_2\in [-\theta_0,\theta_0].
\end{cases}\ee
The last two conditions follow from \eqref{shock20} where $q_i^{\pm} (i=2,3)$ are evaluated at $(\hat{{\bf V}},\hat{V}_7)$ and the compatibility condition \eqref{class2}.

The unique solvability of the problem \eqref{den32} with \eqref{den321} can be verified by several steps using the Duhamel's principle as follows.

{\bf Step 4.1} First, taking the divergence operator for the second, third, and fourth equations in \eqref{den32} leads to
\be\label{den33}\q\q\q\q\q\begin{cases}
(\p_{y_1}^2 +\frac{1}{y_1}\p_{y_1}+ \frac{1}{y_1^2}\p_{y_2}^2+ \p_{y_3}^2) \Pi\\
= \p_{y_1} G_1 +\frac{G_1}{y_1} + \frac{1}{y_1}\p_{y_2} G_2 + \p_{y_3}G_3,&\text{in }\mathbb{D}\\
\Pi(r_s,y')=\Pi(r_2,y')=0,\ &\forall y'\in E,\\
\Pi(y_1,\pm \theta_0,y_3)=0,\ &\text{on }\Sigma_2^{\pm},\\
\Pi(y_1,y_2,\pm 1)= 0\ &\text{on }\Sigma_3^{\pm}.
\end{cases}\ee

There exists a unique $C^{2,\alpha}(\overline{\mathbb{D}})$ smooth solution $\Pi$ to \eqref{den33} with the estimate
\be\no
&&\|\Pi\|_{C^{2,\alpha}(\overline{\mathbb{D}})}\leq C_*\sum_{j=1}^3\|G_j\|_{C^{1,\alpha}(\overline{\mathbb{D}})}
\leq C_*(\epsilon\|(\hat{{\bf V}}, \hat{V}_7)\|_{\Xi}+\|(\hat{{\bf V}}, \hat{V}_7)\|_{\Xi}^2)\leq C_*(\epsilon \delta_0 +\delta_0^2)\no.
\ee
Furthermore, the following compatibility conditions hold
\be\label{cp331}\begin{cases}
\p_{y_1}\Pi(y_1,\pm\theta_0,y_3)=\p_{y_3}\Pi(y_1,\pm\theta_0,y_3)=0,\ \ &\text{on }\Sigma_2^{\pm},\\
\p_{y_1}\Pi(y_1,y_2,\pm 1)=\p_{y_2}\Pi(y_1,y_2,\pm 1)=0,\ \ &\text{on }\Sigma_3^{\pm}.
\end{cases}\ee

{\bf Step 4.2} Next we are going to solve the following divergence-curl system with homogeneous normal boundary conditions
\be\label{den34}\begin{cases}
(\p_{y_1}+\frac1{y_1})\{d\dot{V}_1\}+\frac{1}{y_1}\p_{y_2} \{d_0\dot{V}_2\} + \p_{y_3}\{d_0\dot{V}_3\}=0,\ \ &\text{in }\mathbb{D},\\
\frac{1}{y_1}\p_{y_2} \{d_0\dot{V}_3\}- \p_{y_3} \{d_0\dot{V}_2\}=G_1(\hat{{\bf V}},\hat{V}_7)-\p_{y_1}\Pi:=\tilde{G}_1,\ \ &\text{in }\mathbb{D},\\
\p_{y_3} \{d\dot{V}_1\}-\p_{y_1} \{d_0\dot{V}_3\}= G_2(V_5,V_6;\hat{{\bf V}},\hat{V}_7)-\frac{1}{y_1}\p_{y_2}\Pi:=\tilde{G}_2,\ \ &\text{in }\mathbb{D},\\
(\p_{y_1}+\frac1{y_1})\{d_0\dot{V}_2\}-\frac{1}{y_1}\p_{y_2} \{d\dot{V}_1\}= G_3(V_5,V_6;\hat{{\bf V}},\hat{V}_7)-\p_{y_3}\Pi:=\tilde{G}_3,\ \ &\text{in }\mathbb{D},\\
\dot{V}_1(r_s,y')= \dot{V}_1(r_2,y')=0,\ \ &\forall y'\in E,\\
\dot{V}_2(y_1,\pm\theta_0,y_3)=0,\ \ &\text{on }\Sigma_2^{\pm},\\
\dot{V}_3(y_1,y_2,\pm 1)=0,\ \ &\text{on }\Sigma_3^{\pm}.
\end{cases}\ee

Since $\Pi$ satisfies the equation in \eqref{den33}, then
\be\no
\p_{y_1} \tilde{G}_1+\frac{1}{y_1} \tilde{G}_1+ \frac{1}{y_1}\p_{y_2} \tilde{G}_2+\p_{y_3} \tilde{G}_3\equiv 0,\ \text{in }\mathbb{D}.
\ee
Also it follows from \eqref{cp20} and \eqref{cp331} that
\be\label{cp340}\begin{cases}
(\tilde{G}_1,\tilde{G}_3,\p_{y_2}\tilde{G}_2)(y_1,\pm\theta_0,y_3)=0,\ \ &\text{on }\Sigma_2^{\pm},\\
(\tilde{G}_1,\tilde{G}_2,\p_{y_3}\tilde{G}_3)(y_1,y_2,\pm 1)=0,\ \ &\text{on }\Sigma_3^{\pm}.
\end{cases}
\ee

The unique solvability of the divergence-curl system with the homogeneous normal boundary condition is well-known (cf. \cite{KY2009} and the references therein). By the compatibility condition \eqref{cp340} and the symmetric extension technique as above, there exists a unique $C^{2,\alpha}(\overline{\mathbb{D}})$ smooth vector field $(\dot{V}_1,\dot{V}_2,\dot{V}_3)$ solving \eqref{den34} with
\be\no
&&\sum_{j=1}^3\|\dot{V}_j\|_{C^{2,\alpha}(\overline{\mathbb{D}})}\leq C_*\sum_{j=1}^3\|\tilde{G}_j\|_{C^{1,\alpha}(\overline{\mathbb{D}})}\leq C_*\sum_{j=1}^3\|G_j\|_{C^{1,\alpha}(\overline{\mathbb{D}})}+
\|\Pi\|_{C^{2,\alpha}(\overline{D})}\\\no
&&\leq C_*\sum_{j=1}^3\|G_j\|_{C^{1,\alpha}(\overline{\mathbb{D}})}\leq C_*(\epsilon\|(\hat{{\bf V}}, \hat{V}_7)\|_{\Xi}+\|(\hat{{\bf V}}, \hat{V}_7)\|_{\Xi}^2)\leq C_*(\epsilon \delta_0 +\delta_0^2),
\ee
and the following compatibility conditions hold
\be\label{den342}\begin{cases}
\p_{y_2} (\dot{V}_1,\dot{V}_3)(y_1,\pm\theta_0,y_3)=(\dot{V}_2,\p_{y_2}^2 \dot{V}_2)(y_1,\pm\theta_0,y_3)=0, \ &\text{on }\Sigma_2^{\pm},\\
\p_{y_3} (\dot{V}_1,\dot{V}_2)(y_1,y_2,\pm 1)=(\dot{V}_3,\p_{y_3}^2 \dot{V}_3)(y_1,y_2, \pm 1)=0, \ &\text{on }\Sigma_3^{\pm}.
\end{cases}\ee

{\bf Step 4.3} Let $(V_1,V_2,V_3)$ be the solution to \eqref{den32}, and set
\be\no
N_j(y)= V_j(y)- \dot{V}_j(y), j=1,2,3.
\ee
Then $N_j, j=1,2,3$ solve the following equations
\be\label{den36}\begin{cases}
d_1\p_{y_1} N_1+\frac{1}{y_1}\p_{y_2} N_2 +\p_{y_3} N_3 +\frac{1}{y_1} N_1+d_2 N_1=G_4(\hat{{\bf V}},\hat{V}_7),\\
\frac{1}{y_1}\p_{y_2} \{d_0N_3\}-\p_{y_3} \{d_0 N_2\}=0,\\
\p_{y_3} \{d N_1+\frac{a_2}{a_1}d_3(y_1) N_1(r_s,y')\}-\p_{y_1} \{d_0 N_3\}= 0,\\
(\p_{y_1}+\frac{1}{y_1})\{d_0 N_2\}-\frac{1}{y_1}\p_{y_2} \{d N_1+\frac{a_2}{a_1}d_3(y_1) N_1(r_s,y')\}= 0,\\
\end{cases}\ee
endowed with the boundary conditions:
\be\label{den360}\begin{cases}
N_2(y_1,\pm \theta_0,y_3)=0,\ &\text{on }\Sigma_2^{\pm},\\
N_3(y_1,y_2,\pm 1)=0,\ &\text{on }\Sigma_3^{\pm},\\
d(r_2)N_1(r_2,y')+\frac{a_2}{a_1} d_{3}(r_2) N_1(r_s,y')=q_4(\hat{{\bf V}}(r_s,y'),\hat{V}_7(y')),  &\forall y'\in E,\\
(\frac{1}{r_s^2}\p_{y_2}^2+\p_{y_3}^2) N_1(r_s,y')-a_0 a_1 (\frac{1}{r_s}\p_{y_2} N_2+\p_{y_3}N_3)(r_s,y')=q_5(y'), &\forall y'\in E,\\
(\frac{1}{r_s}\p_{y_2}N_1-a_0a_1N_2)(r_s,\pm \theta_0,y_3)=0,\ &\forall y_3\in [-1,1],\\
(\p_{y_3} N_1-a_0 a_1 N_3)(r_s,y_2,\pm 1)= 0,\ &\forall y_2\in [-\theta_0,\theta_0].
\end{cases}\ee
where
\be\no
&&G_4(\hat{{\bf V}},\hat{V}_7)=d_5 V_5+G_0(\hat{{\bf V}},\hat{V}_7)+
\left(\frac{d(y_1)}{d_0(y_1)}-d_1(y_1)\right)\p_{y_1}\dot{V}_1\\
&&\quad\quad\quad\quad\,\, +\left( \frac{d(y_1)}{d_0(y_1) y_1}-\frac{1}{y_1}+\dfrac{d'(y_1)}{d_0(y_1)}-d_2(y_1) \right)\dot{V}_1
,\no\\
&&q_5(\hat{{\bf V}}(r_s,y'),\hat{V}_7(y'))=q_1(\hat{{\bf V}}(r_s,y'),\hat{V}_7(y'))+a_0 a_1(\frac{1}{r_s}\p_{y_2}\dot{V}_2+\p_{y_3} \dot{V}_3)(r_s,y').\no
\ee

It follows from the second, third, and fourth equations in \eqref{den36} that there exists a potential function $\phi$ such that
\be\no
&&d(y_1)N_1(y_1,y') +\frac{a_2}{a_1} d_3(y_1) N_1(r_s,y')=\p_{y_1}\phi(y_1,y'),\\\no
&&d_0(y_1)N_2(y_1,y')=\frac{1}{y_1}\p_{y_2}\phi(y_1,y'), \ d_0(y_1)N_3=\p_{y_3}\phi(y_1,y').
\ee
Therefore
\be\no
&&N_1(r_s,y')=\frac 1{a_3}\p_{y_1}\phi(r_s,y'),\\\no
&&N_1(y_1,y')=\frac{1}{d(y_1)}\p_{y_1}\phi(y_1,y')- \frac{a_2}{a_1a_3} \frac{d_3}{d}(y_1)\p_{y_1}\phi(r_s,y'),\\\no
&&N_2(y_1,y')=\frac{1}{d_0(y_1)}\frac{1}{y_1}\p_{y_2}\phi,\ \ N_3(y_1,y')=\frac{1}{d_0(y_1)}\p_{y_3}\phi,
\ee
with
\be\no
a_3=\frac{(\gamma-1) \bar{M}^2(r_s)+1}{\gamma \bar{M}^2(r_s)}>0.
\ee
Thus the problem \eqref{den36} and \eqref{den360} is equivalent to
\be\label{den37}\q\begin{cases}
d_1(y_1)\p_{y_1}(\frac{\p_{y_1}\phi}{d(y_1)}) +\frac{1}{d_0(y_1)}(\frac{1}{y_1^2}\p_{y_2}^2 \phi +\p_{y_3}^2 \phi) +(\frac{1}{y_1}+ d_2(y_1))\frac{\p_{y_1}\phi}{d(y_1)}\\
\q\q-\frac{a_2}{a_1 a_3} d_4(y_1)\p_{y_1}\phi(r_s,y')=G_4(\hat{{\bf V}},\hat{V}_7),\ \ &\text{in }\mathbb{D},\\
\p_{y_2}\phi(y_1,\pm \theta_0,y_3)=0,\ \ &\text{on }\Sigma_2^{\pm},\\
\p_{y_3}\phi(y_1,y_2,\pm 1)= 0,\ \ \ \ &\text{on }\Sigma_3^{\pm},\\
\p_{y_1}\phi(r_2,y')=q_4(y'),\ \ &\forall y'\in E,\\
(\frac{1}{r_s^2}\p_{y_2}^2+\p_{y_3}^2)\big(\p_{y_1}\phi(r_s,y')-a_4 \phi(r_s,y')\big)=a_3q_5(y'),\ \ &\forall y'\in E,\\
\frac{1}{r_s}\p_{y_2}(\p_{y_1}\phi-a_4 \phi)(r_s,\pm \theta_0,y_3)=0,\ &\forall y_3\in [-1,1],\\
\p_{y_3} (\p_{y_1}\phi-a_4\phi)(r_s,y_2,\pm 1)=0,\ \ &\forall y_2\in [-\theta_0,\theta_0].
\end{cases}\ee
where
\be\no
&&d_4(y_1)=d_1(y_1)\frac{d}{d y_1}\bigg(\frac{d_3(y_1)}{d(y_1)}\bigg)+(\frac{1}{y_1}+d_2(y_1)) \frac{d_3(y_1)}{d(y_1)},\no\\
&& a_4=a_0 a_1 a_3>0.\no
\ee

A direct computation shows that
\begin{equation*}
\begin{split}
d'(y_1)&=\frac{3+(\ga-2)\bM_+^2}{y_1(\bM_+^2-1)}\brho_+ \bka^2 \bM_+^2=-\brho_+\bka^2d_2 -\frac{1}{y_1} \brho_+ \bka^2 \bM_+^2,\\
d_3'(y_1)&=\frac{2\bB+\bU_+^2}{2y_1(1-\bM_+^2)\ga S_+ \bU_+}-\frac{\bka^2\brho_+ \bU_+}{y_1\bS_+(\ga-1)}\\
&=\frac{d_3}{y_1d_1}-\frac{\bka^2\brho_+ \bU_+}{(\ga-1)\bS_+}\frac{1}{y_1}(1+\frac{1}{d_1})+\frac{1}{y_1d_1}\frac{\bU_+}{\ga \bS_+}.
\end{split}
\end{equation*}

As a result, we claim that $d_4(y_1)>0$ for any $y_1\in[r_s,r_2]$. Indeed, it holds that
\begin{equation}
\begin{split}
d_4(y_1)&=d_1(y_1)(\frac{d_3(y_1)}{d(y_1)})'+(\frac{1}{y_1}+d_2(y_1))\frac{d_3(y_1)}{d(y_1)}\\
&=\frac{d_1}{d^2}(d_3'd-d_3 d')+\frac{d_3}{y_1 d}+\frac{d_2 d_3}{d}\\
&=\frac{d_3}{y_1d}-\frac{\bka^2\brho_+ \bU_+}{(\ga-1)\bS_+}\frac{1}{y_1d}(1+d_1)+\frac{1}{y_1d}\frac{\bU_+}{\ga \bS_+}\\
&+\frac{d_1 d_3}{d^2}(\brho_+\bka^2\,d_2 +\frac{1}{y_1} \brho_+ \bka^2 \bM_+^2)+\frac{d_3}{y_1 d}+\frac{d_2 d_3}{d}\\
&=\frac{1}{y_1d}\left(
\frac{2\bB}{\ga \bS_+ \bU_+}+\frac{\bka^2 \brho_+ \bU_+}{(\ga-1)\bS_+}\bM_+^2\right)+\frac{d_2 d_3}{d}\\
&+\frac{d_1 d_3}{d^2}(\brho_+\bka^2\,d_2 +\frac{1}{y_1} \brho_+ \bka^2 \bM_+^2)>0,
\end{split}
\end{equation}
because $d(y_1)>0$ and $d_i(y_1)>0,\,i=1,2,3$, for all $y_1\in[r_s,r_2]$.

Resolving the Poisson equation with the Neumann boundary conditions in the last three equations in \eqref{den37}, we derive an oblique boundary condition for the potential $\phi$ on the boundary $\{(r_s,y'): y'\in E\}$.
\begin{lemma}\label{oblique}({\bf The oblique boundary condition on the shock front.})
On the shock front $\{(r_s, y'): y'\in E\}$, there exists a unique $C^{2,\alpha}(\overline{E})$ function $m_1(y')$ such that
\be\no
\p_{y_1}\phi(r_s,y')-a_4 \phi(r_s,y')= m_1(y'),
\ee
where $m_1(y')$ satisfies the Poisson equation with the homogeneous Neumann boundary conditions
\be\label{den38}\begin{cases}
(\frac{1}{r_s^2}\p_{y_2}^2+\p_{y_3}^2)m_1(y')=a_3q_5(\hat{{\bf V}}(r_s,y'),\hat{V}_7(y')),\ \ &\text{in } E,\\
\frac{1}{r_s}\p_{y_2}m_1(\pm \theta_0,y_3)=0,\ &\forall y_3\in [-1,1],\\
\p_{y_3} m_1(y_2,\pm 1)=0,\ \ &\forall y_2\in [-\theta_0,\theta_0],
\end{cases}\ee
and the condition
\be\label{den39}
\iint_{E} m_1(y') dy'=0.
\ee
\end{lemma}

Then the problem \eqref{den37} can be reduced to
\be\label{den41}\q\begin{cases}
d_1(y_1)\p_{y_1}(\frac{\p_{y_1}\phi}{d(y_1)}) +\frac{1}{d_0(y_1)}(\frac{1}{y_1^2}\p_{y_2}^2 \phi +\p_{y_3}^2 \phi)\\\q\q\q\q +(\frac{1}{y_1}+ d_2(y_1))\frac{\p_{y_1}\phi}{d(y_1)}-a_0a_2 d_4\phi(r_s,y')=G_5(y),\ &\text{in }\mathbb{D},\\
\p_{y_1}\phi(r_s,y')-a_4 \phi(r_s,y')=m_1(y'),  &\forall y'\in E,\\
\p_{y_1}\phi(r_2,y')=m_2(y'),  &\forall y'\in E,\\
\p_{y_2}\phi(y_1,\pm \theta_0,y_3)=0,   &\text{on }\Sigma_2^{\pm},\\
\p_{y_3}\phi(y_1,y_2,\pm 1)= 0,  &\text{on }\Sigma_3^{\pm},
\end{cases}\ee
where
\be\no
&&G_5(y)=G_4(y)+ \frac{a_2}{a_1 a_3} d_4(y_1) m_1(y'),\\\no
&&m_2(y')= q_4(\hat{{\bf V}}(r_s,y'),\hat{V}_7(y')).
\ee
It can be checked easily that the function $G_5$ and $m_i, i=1,2$ satisfy the following compatibility conditions
\be\label{den415}\begin{cases}
\p_{y_2}G_5(y_1,\pm\theta_0, y_3)=0,\ &\text{on }\Sigma_2^{\pm},\\
\p_{y_3}G_5(y_1,y_2, \pm 1)=0,\ &\text{on }\Sigma_3^{\pm},\\
\p_{y_2} m_1(\pm \theta_0,y_3)=\p_{y_2} m_2(\pm \theta_0,y_3)=0, \ &\forall y_3\in [-1,1],\\
\p_{y_3} m_1(y_2,\pm 1)=\p_{y_3} m_2(y_2, \pm 1)=0, \ &\forall y_2\in [-\theta_0,\theta_0].
\end{cases}\ee

Then under the assumptions \eqref{supA2}, we have
\be\label{dsign}\begin{cases}
d_{1}(y_1)>0, d_{0}(y_1)=\bar{\kappa}^2\bar{\rho}_+(y_1)(\bar{A}_+^2(y_1)-1)>0,\ \ &\forall y_1\in [r_s,r_2],\\
d(y_1)=\bar{\kappa}^2\bar{\rho}_+(y_1)(\bar{A}_+^2(y_1)+\bar{M}_+^2(y_1)-1)>0,\ \ &\forall y_1\in [r_s,r_2].
\end{cases}\ee
Thus the equation in \eqref{den41} is a second order elliptic equation with a nonlocal term. The problem \eqref{den41} can be solved similarly to the method used in \cite[Proposition 3.4]{WengX24ArxicCyl}.

\begin{proposition}\label{solvability}
{\it Suppose that $G_5\in C^{1,\alpha}(\overline{\mathbb{D}})$ and $(m_1,m_2)\in (C^{2,\alpha}(\overline{E}))^2$ satisfy \eqref{den415}. Then there exists a unique $C^{3,\alpha}(\overline{\mathbb{D}})$ solution to the problem \eqref{den41} with
\be\label{den414}
\|\phi\|_{C^{3,\alpha}(\overline{\mathbb{D}})}\leq C_*(\|G_5\|_{C^{1,\alpha}(\overline{\mathbb{D}})}+\sum_{j=1}^2\|m_j\|_{C^{2,\alpha}(\overline{E})}),
\ee
where $C_*$ depends only on $d_1,d_4,d_5, a_3,a_4$ and thus on the background solution.
}\end{proposition}

Thus $N_1(y)=\p_{y_1}\phi(y)-\frac{1}{a_3} d_3(y_1)\p_{y_1}\phi(r_s,y'), N_2(y)=\frac{1}{y_1}\p_{y_2}\phi(y_1,y')$ and $N_3(y)=\p_{y_3}\phi$ would solve the problem \eqref{den36}-\eqref{den360}. Differentiating the first equation in \eqref{den36} with respect to $y_2$ (resp. $y_3$) and evaluating at $y_2=\pm\theta_0$ (resp. $y_3=\pm 1$), one gets from \eqref{den342} that
\be\no\begin{cases}
\p_{y_2}^2 N_2(y_1,\pm \theta_0,y_3)=0,\ \ &\text{on } \Sigma_2^{\pm},\\
\p_{y_3}^2 N_3(y_1,y_2,\pm 1)=0,\ \ &\text{on } \Sigma_3^{\pm}.
\end{cases}\ee

Then
\be\no
&&V_1(y_1,y')=\dot{V}_1(y)+\p_{y_1}\phi(y)-\frac{1}{a_3} d_3(y_1)\p_{y_1}\phi(r_s,y'),\\\no
&&V_2(y_1,y')=\dot{V}_2(y_1,y')+\frac{1}{y_1}\p_{y_2}\phi(y_1,y'),\no\\
&&V_3(y_1,y')=\dot{V}_3(y_1,y')+\p_{y_3}\phi(y_1,y')\no,
\ee
will solve the problem \eqref{den32} with \eqref{den321} and satisfy
\be\no
&&\sum_{j=1}^3\|V_j\|_{C^{2,\alpha}(\overline{\mathbb{D}})}\leq C_*(\sum_{j=1}^3\|\dot{V}_j\|_{C^{2,\alpha}(\overline{\mathbb{D}})}+\|\nabla \phi\|_{C^{2,\alpha}(\overline{\mathbb{D}})}+ \|\p_{y_1}\phi(r_s,y')\|_{C^{2,\alpha}(\overline{E})})\\\label{den46}
&&\leq C_*(\epsilon +C_*(\epsilon\|(\hat{{\bf V}}, \hat{V}_7)\|_{\Xi}+\|(\hat{{\bf V}}, \hat{V}_7)\|_{\Xi}^2)\leq C_*(\epsilon+\epsilon \delta_0 +\delta_0^2).
\ee
Also the following compatibility conditions hold
\be\label{den461}\begin{cases}
(V_2,\p_{y_2}^2 V_2,\p_{y_2} V_1,\p_{y_2} V_3)(y_1,\pm\theta_0,y_3)=0,\ \ \ &\text{on } \Sigma_2^{\pm},\\
(V_3,\p_{y_3}^2 V_3,\p_{y_3} V_1,\p_{y_3} V_2)(y_1,y_2,\pm 1)=0,\ \ &\text{on } \Sigma_3^{\pm}.
\end{cases}\ee

{\bf Step 5.} Once the velocity fields $V_1$, $V_2$, and $V_3$ are obtained, the function $V_4$ in \eqref{ent43} can be uniquely determined as follows.
\be\label{ent53}
V_4(y_1,y')=\frac{a_2}{a_1} V_1(r_s,y')+R_4(\hat{{\bf V}}(r_s,\beta_2(y),\beta_3(y)),\hat{V}_7(\beta_2(y),\beta_3(y))).
\ee
Then it can be checked easily that the following estimate and compatibility conditions hold:
\be\no
&&\|V_4\|_{C^{2,\alpha}(\overline{\mathbb{D}})}\leq C_*\|V_1(r_s,\cdot)\|_{C^{2,\alpha}(\overline{E})}+ C_*(\epsilon \|(\hat{{\bf V}}, \hat{V}_7)\|_{\Xi}+\|(\hat{{\bf V}}, \hat{V}_7)\|_{\Xi}^2)\\\label{ent54}
&&\leq C_*(\epsilon \delta_0+ \delta_0^2),
\ee
and
\be\label{ent55}\begin{cases}
\p_{y_2} V_4(y_1,\pm \theta_0, y_3)=\frac{a_2}{a_1}\p_{y_2} V_1(r_s,\pm\theta_0,y_3)=0,\ \ &\text{on }\Sigma_2^{\pm},\\
\p_{y_3} V_4(y_1,y_2,\pm 1)=\frac{a_2}{a_1}\p_{y_3} V_1(r_s,y_2,\pm 1)=0,\ \ &\text{on }\Sigma_3^{\pm}.
\end{cases}\ee

Finally, the shock front is given by
\be\label{shock50}
V_7(y')=\frac{1}{a_1} V_1(r_s,y')-\frac{1}{a_1} R_1(\hat{{\bf V}}(r_s,y'),\hat{V}_7(y')),
\ee
and it is clear that $V_7\in C^{2,\alpha}(\overline{E})$ and
\be\label{shock51}\begin{cases}
\p_{y_2}V_7(\pm\theta_0, y_3)=0,\ \ &\text{on }y_3\in [-1,1],\\
\p_{y_3}V_7(y_2, \pm 1)=0,\ \ &\text{on }y_2\in [-\theta_0,\theta_0].
\end{cases}\ee

Furthermore, there holds
\be\label{shock53}\begin{cases}
\frac{1}{r_s}\p_{y_2} V_7 (y')= \frac{a_0}{a_1} V_2(r_s,y') + g_2(\hat{{\bf V}}(r_s,y'),\hat{V}_7(y')), \ \ &\text{in }E,\\
\p_{y_3} V_7 (y')= \frac{a_0}{a_1} V_3(r_s,y') + g_3(\hat{{\bf V}}(r_s,y'),\hat{V}_7(y')), \ \ &\text{in }E.
\end{cases}\ee
Therefore $V_7\in C^{3,\alpha}(\overline{E})$ admits the following estimate
\be\label{shock54}
&&\|V_7\|_{C^{3,\alpha}(\overline{E})}\leq C_*(\|V_1(r_s,\cdot)\|_{C^{2,\alpha}(\overline{E})}+ \|R_1(\hat{{\bf V}}(r_s,y'),\hat{V}_7(y'))\|_{C^{2,\alpha}(\overline{E})})\\\no
&&\quad + C_*\sum_{j=2}^3(\|V_j(r_s,\cdot)\|_{C^{2,\alpha}(\overline{E})}+\|g_j(\hat{{\bf V}}(r_s,y'),\hat{V}_7(y'))\|_{C^{2,\alpha}(\overline{E})})\\\no
&&\leq C_*(\epsilon + \epsilon  \|(\hat{{\bf V}}, \hat{V}_7)\|_{\Xi}+\|(\hat{{\bf V}}, \hat{V}_7)\|_{\Xi}^2)\leq C_*(\epsilon +\epsilon \delta_0+\delta_0^2),
\ee
and
\be\label{shock55}\begin{cases}
\p_{y_2}^3V_7(\pm\theta_0, y_3)=0,\ \ &\forall y_3\in [-1,1],\\
\p_{y_3}^3 V_7(y_2, \pm 1)=0,\ \ &\forall y_2\in [-\theta_0,\theta_0].
\end{cases}\ee

Combining the estimates \eqref{ber43}, \eqref{den46}, \eqref{ent54} and \eqref{shock54}, one concludes that
\be\no
\|({\bf V}, V_7)\|_{\Xi}= \sum_{j=1}^6 \|V_j\|_{C^{2,\alpha}(\overline{\mathbb{D}})}+\|V_7\|_{C^{3,\alpha}(\overline{E})}\leq C_*(\epsilon +\epsilon \delta_0+\delta_0^2)\leq C_*(\epsilon +\delta_0^2).
\ee
Choose $\delta_0=\sqrt{\epsilon}$ and let $\epsilon<\epsilon_0=\frac{1}{4 C_*^2}$. Then $\|({\bf V}, V_7)\|_{\Xi}\leq 2C_*\epsilon\leq \delta_0$. Furthermore, the compatibility conditions \eqref{ber44}, \eqref{den461},\eqref{ent55}, \eqref{shock51} and \eqref{shock55} hold, thus $({\bf V}, V_7)\in \Xi$. We now can define the operator $\mathcal{T}:(\hat{{\bf V}}, \hat{V}_7)\mapsto ({\bf V}, V_7)$ which maps $\Xi$ to itself.

{\bf Step 6.} The contraction of the operator $\mathcal{T}$ can be proved similarly as in Step 6 \cite[Theorem 2.3]{WengX24ArxicCyl}. Thus  $\mathcal{T}$ has a unique fixed point $({\bf V}, V_7) \in \Xi$. Furthermore, it can be proven that the auxiliary function $\Pi$, introduced in equations \eqref{den32}--\eqref{den321}, also vanishes for the fixed point $({\bf V}, V_7)$. The proof is completed. 
%
%

\section*{Acknowledgment}
Weng is supported by National Natural Science Foundation of China 12571240, 12221001.

{\bf Data Availability Statement.} No data, models or code were generated or used during the study.

{\bf Conflict of interest.} On behalf of all authors, the corresponding author states that there is no conflict of interests.

\normalem
\bibliographystyle{siam}
\bibliography{MHDShocks}

\begin{thebibliography}{10}

\bibitem{BE1959}
{\sc J.~Bazer and W.~B. Ericson}, {\em Hydromagnetic shocks}, The Astrophysical
  Journal, 129 (1959), p.~758.

\bibitem{BS2007}
{\sc S.~Benzoni-Gavage and D.~Serre}, {\em Multidimensional hyperbolic partial
  differential equations}, Oxford Mathematical Monographs, The Clarendon Press,
  Oxford University Press, Oxford, 2007.
\newblock First-order systems and applications.

\bibitem{BT2002}
{\sc A.~Blokhin and Y.~Trakhinin}, {\em Stability of strong discontinuities in
  fluids and {MHD}}, in Handbook of mathematical fluid dynamics, {V}ol. {I},
  North-Holland, Amsterdam, 2002, pp.~545--652.

\bibitem{ChenF2003}
{\sc G.-Q. Chen and M.~Feldman}, {\em Multidimensional transonic shocks and
  free boundary problems for nonlinear equations of mixed type}, J. Amer. Math.
  Soc., 16 (2003), pp.~461--494.

\bibitem{ChenW2008}
{\sc G.-Q. Chen and Y.-G. Wang}, {\em Existence and stability of compressible
  current-vortex sheets in three-dimensional magnetohydrodynamics}, Arch.
  Ration. Mech. Anal., 187 (2008), pp.~369--408.

\bibitem{Chen2008}
{\sc S.~Chen}, {\em Transonic shocks in 3-{D} compressible flow passing a duct
  with a general section for {E}uler systems}, Trans. Amer. Math. Soc., 360
  (2008), pp.~5265--5289.

\bibitem{ChenYuan2008}
{\sc S.~Chen and H.~Yuan}, {\em Transonic shocks in compressible flow passing a
  duct for three-dimensional {E}uler systems}, Arch. Ration. Mech. Anal., 187
  (2008), pp.~523--556.

\bibitem{Chu1962}
{\sc C.~K. Chu}, {\em Magnetohydrodynamic nozzle flow with three transitions},
  The Physics of Fluids, 5 (1962), pp.~550--559.

\bibitem{CF1948}
{\sc R.~Courant and K.~O. Friedrichs}, {\em Supersonic {F}low and {S}hock
  {W}aves}, Interscience Publishers, Inc., New York, 1948.

\bibitem{Davidson2017}
{\sc P.~A. Davidson}, {\em Introduction to magnetohydrodynamics}, Cambridge
  Texts in Applied Mathematics, Cambridge University Press, Cambridge,
  second~ed., 2017.

\bibitem{DT1950}
{\sc F.~De~Hoffmann and E.~Teller}, {\em Magnetohydrodynamic shocks}, Physical
  Review, 80 (1950), pp.~692--703.

\bibitem{FangX2021}
{\sc B.~Fang and Z.~Xin}, {\em On admissible locations of transonic shock
  fronts for steady {E}uler flows in an almost flat finite nozzle with
  prescribed receiver pressure}, Comm. Pure Appl. Math., 74 (2021),
  pp.~1493--1544.

\bibitem{FT2013}
{\sc G.~D. Fleishman and I.~N. Toptygin}, {\em Cosmic Electrodynamics:
  Electrodynamics and Magnetic Hydrodynamics of Cosmic Plasmas}, Astrophysics
  and Space Science Library, Springer New York, NY, 2013.

\bibitem{Fried1954}
{\sc K.~O. Friedrichs}, {\em Nonlinear wave motion in magnetohydrodynamics},
  Technical Report LAMS-2105, Los Alamos Scientific Laboratory, Los Alamos, NM,
  USA, 1954.

\bibitem{Grad1960}
{\sc H.~Grad}, {\em Reducible problems in magneto-fluid dynamic steady flows},
  Reviews of Modern Physics, 32 (1960), pp.~830--847.

\bibitem{KY2009}
{\sc H.~Kozono and T.~Yanagisawa}, {\em {$L^r$}-variational inequality for
  vector fields and the {H}elmholtz-{W}eyl decomposition in bounded domains},
  Indiana Univ. Math. J., 58 (2009), pp.~1853--1920.

\bibitem{LXY2009CMP}
{\sc J.~Li, Z.~Xin, and H.~Yin}, {\em On transonic shocks in a nozzle with
  variable end pressures}, Comm. Math. Phys., 291 (2009), pp.~111--150.

\bibitem{LXY2010JDE}
\leavevmode\vrule height 2pt depth -1.6pt width 23pt, {\em On transonic shocks
  in a conic divergent nozzle with axi-symmetric exit pressures}, J.
  Differential Equations, 248 (2010), pp.~423--469.

\bibitem{LXY2013}
\leavevmode\vrule height 2pt depth -1.6pt width 23pt, {\em Transonic shocks for
  the full compressible {E}uler system in a general two-dimensional de {L}aval
  nozzle}, Arch. Ration. Mech. Anal., 207 (2013), pp.~533--581.

\bibitem{LiuXY2016}
{\sc L.~Liu, G.~Xu, and H.~Yuan}, {\em Stability of spherically symmetric
  subsonic flows and transonic shocks under multidimensional perturbations},
  Adv. Math., 291 (2016), pp.~696--757.

\bibitem{MTT2018}
{\sc A.~Morando, Y.~Trakhinin, and P.~Trebeschi}, {\em Local existence of {MHD}
  contact discontinuities}, Arch. Ration. Mech. Anal., 228 (2018),
  pp.~691--742.

\bibitem{Trak2003CMP}
{\sc Y.~Trakhinin}, {\em A complete 2{D} stability analysis of fast {MHD}
  shocks in an ideal gas}, Comm. Math. Phys., 236 (2003), pp.~65--92.

\bibitem{Tark2009ARMA}
\leavevmode\vrule height 2pt depth -1.6pt width 23pt, {\em The existence of
  current-vortex sheets in ideal compressible magnetohydrodynamics}, Arch.
  Ration. Mech. Anal., 191 (2009), pp.~245--310.

\bibitem{TrakWang2022}
{\sc Y.~Trakhinin and T.~Wang}, {\em Nonlinear stability of {MHD} contact
  discontinuities with surface tension}, Arch. Ration. Mech. Anal., 243 (2022),
  pp.~1091--1149.

\bibitem{WangXin2024}
{\sc Y.~Wang and Z.~Xin}, {\em Existence of multi-dimensional contact
  discontinuities for the ideal compressible magnetohydrodynamics}, Comm. Pure
  Appl. Math., 77 (2024), pp.~583--629.

\bibitem{Weng25ArxivSph}
{\sc S.~Weng}, {\em Three dimensional spherical transonic shock in a
  hemispherical shell},  (2025).
\newblock arXiv:2503.14886.

\bibitem{WengXieXin2021}
{\sc S.~Weng, C.~Xie, and Z.~Xin}, {\em Structural stability of the transonic
  shock problem in a divergent three-dimensional axisymmetric perturbed
  nozzle}, SIAM J. Math. Anal., 53 (2021), pp.~279--308.

\bibitem{WengX2019}
{\sc S.~Weng and Z.~Xin}, {\em A deformation-curl decomposition for three
  dimensional steady euler equations}, Scientia Sinica Mathematica, 49 (2019),
  pp.~307--320.
\newblock (in Chinese).

\bibitem{WengX24ArxicCyl}
{\sc S.~Weng and Z.~Xin}, {\em Existence and stability of cylindrical transonic
  shock solutions under three dimensional perturbations},  (2024).
\newblock arXiv:2304.02429.

\bibitem{WengY2024}
{\sc S.~Weng and W.~Yang}, {\em Structural stability of transonic shock flows
  with an external force}, Proc. Roy. Soc. Edinburgh,  (2024), p.~1–24.

\bibitem{WengZZ2025}
{\sc S.~Weng, Z.~Zhang, and Y.~Zhou}, {\em Structural stability of three
  dimensional transonic shock flows with an external force}, J. Differential
  Equations, 427 (2025), pp.~310--349.

\bibitem{XinYin2005}
{\sc Z.~Xin and H.~Yin}, {\em Transonic shock in a nozzle. {I}.
  {T}wo-dimensional case}, Comm. Pure Appl. Math., 58 (2005), pp.~999--1050.

\bibitem{XinYin2008PJM}
\leavevmode\vrule height 2pt depth -1.6pt width 23pt, {\em Three-dimensional
  transonic shocks in a nozzle}, Pacific J. Math., 236 (2008), pp.~139--193.

\bibitem{XinYin2008JDE}
\leavevmode\vrule height 2pt depth -1.6pt width 23pt, {\em The transonic shock
  in a nozzle, 2-{D} and 3-{D} complete {E}uler systems}, J. Differential
  Equations, 245 (2008), pp.~1014--1085.

\end{thebibliography}

\end{document}